\documentclass[a4paper,11pt,fleqn]{article}
\usepackage{geometry,graphicx,wrapfig,esvect,subfigure,amsmath,amssymb,amsthm,verbatim,sidecap,hyperref} 
\newcommand{\C}{\mathbb{C}}
\newcommand{\Z}{\mathbb{Z}}

\DeclareMathOperator*{\Arg}{Arg}

\title{The Equilateral Pentagon at Zero Angular Momentum: Maximal Rotation Through Optimal Deformation}
\author{William Tong and Holger R. Dullin \\
School of Mathematics and Statistics \\
University of Sydney, Australia}

\begin{document}
\pagenumbering{arabic}
\maketitle

\section*{Abstract}
A pentagon in the plane with fixed side-lengths has a two-dimensional shape space. Considering the pentagon as 
a mechanical system with point masses at the corners we answer the question of how much the 
pentagon can rotate with zero angular momentum. We show that the shape space of the equilateral 
pentagon has genus 4 and find a fundamental region by discrete symmetry reduction with respect to 
symmetry group $D_5$. The amount of rotation $\Delta \theta$ for a loop in shape space at zero angular momentum 
is interpreted as a geometric phase and is obtained as an integral of a function $B$ over the  region of shape space 
enclosed by the loop. With a simple variational argument we determine locally optimal loops as the zero contours 
of the function $B$. The resulting shape change is represented as a Fourier series, and the global 
maximum of $\Delta \theta \approx 45^\circ$ is found for a loop around the regular pentagram. 
We also show that restricting allowed shapes to convex pentagons the optimal loop is the boundary 
of the convex region and gives $\Delta \theta \approx 19^\circ$.

\section{Introduction}
The possibility of achieving overall rotation at zero total angular momentum in an isolated mechanical system is surprising. 
It is possible for non-rigid bodies, in particular for systems of coupled rigid bodies, to change their orientations without an external torque using only internal forces, thus preserving the total angular momentum.
The classical example of this phenomenon is exhibited by falling cats: a cat dropped upside-down without angular momentum will reorient itself by changing its shape and land on its feet, with roughly the same final and initial shape. The first theoretical explanation was given by Kane and Scher 
\cite{KaneScher69}, also see \cite{Montgomery93cat,PuttermanRaz08}.

Here we present a study of the equilateral pentagon in the plane where we permit the angles $\psi_i$, $i = 1, \dots, 5$ 
of adjacent edges to continuously change while the lengths of the edges are all fixed to 1.
The reason the pentagon is chosen among other polygons is because it has a two-dimensional shape space.
By contrast, an equilateral triangle has a fixed shape; an equilateral quadrilateral can change from a square, through a rhombus, 
to a degenerate shape of a line. Adding one more degree of freedom makes it possible for the equilateral pentagon to achieve overall rotation at zero angular momentum through a periodic shape change.
In general the side-lengths of the pentagon could be considered as parameters, but we restrict ourselves to the equilateral case,
which gives some additional simplification and beauty through its discrete symmetry.

Changes in size are irrelevant to our problem, so we use the word {\em shape} in the sense of congruence:
two pentagons have the same shape if one can be transformed into the other by isometries of the Euclidean plane,
that is a combination of rotation, translation and reflections, generating $E(2)$. 
Sometimes we will consider direct isometries $SE(2)$ only, omitting the reflections. 
Thus the equilateral pentagon is a mechanical system with symmetry. 
In this setting the overall rotation at angular momentum zero appears as a {\em geometric phase}.
Symmetry reduction splits the dynamics into a motion in the symmetry group (translations and rotations) and a reduced system. 
The motion in the reduced system drives the motion in the group direction, 
and it is possible (in fact typical) that traversing a closed loop in the reduced system does not lead to a closed orbit in the group,
see for example~\cite{MMR90,LittlejohnReinsch97}. The motion along the group can be split into a geometric and a dynamic phase, where the 
geometric phase does not depend on the speed at which the loop in the reduced system is traversed.
In our case the translation is removed by going to the centre of mass frame, and since we are in the plane
only a single angle $\theta$ is needed to describe the orientation.
Although the inspiration for this work was taken from the general modern theory of geometric phase
\cite{MMR90,LittlejohnReinsch97,Bloch03}, here we take an approach that can be understood with a 
minimal background in mechanics.

The equilateral pentagon as a symmetry reduced mechanical system would perform a certain motion in (the cotangent bundle of) 
shape space if the angles were free to move. 
In many ways the system would then be similar to the 3-linkage studied in \cite{MacKay03}, 
the four-bar linkage studied in \cite{YangKrish},
or, e.g., the planar skater studied in \cite{Marsden98symmetriesin}.
Here, however, we take the point of view that the angles can be completely controlled
by us, for instance with a motor per joint. The only constraint imposed is that the motion must be such that the total angular momentum 
remains constant (at value 0, in particular). 
A way to picture this is to think of the pentagon as a space station that has five motors at the joints and is floating in space 
without angular momentum. 
Controlling the motors, we are free to prescribe any motion in shape space.
We can then ask ``What is the optimal periodic shape change of the equilateral pentagon
so that the overall rotation $\Delta \theta$ after one traversal of the loop in shape space is as large as possible?''
Specifically we seek the global maximum of $\Delta \theta$ on the space of all finite smooth contractible loops.

The plan of the paper is as follows: 
\begin{itemize} \setlength{\itemsep}{0mm}
\item {\em Section 2 -} We describe the equilateral pentagon and its shape space; 
\item {\em Section 3 -} Using reduction by the discrete symmetry group $D_5$, we obtain a fundamental region of shape space and show how all of shape space is tiled by this fundamental piece;
\item {\em Section 4 -}  Explicit formulas for the moment of inertia, angular momentum and the rate of change of orientation are derived;
\item {\em Section 5 -} The geometric phase $\Delta \theta$ is defined in terms of a line integral
which is then converted into an area integral over the enclosed region of a scalar function $B$ on shape space;
\item {\em Section 6 -} We show that the zero level of the function $B$ gives the optimal loop, and we obtain a representation of the corresponding shape change in terms of a Fourier series;
\item {\em Section 7 -}  Restricting to convex pentagons we show that the optimal loop for this sub-family is given by the boundary of the region of convex pentagons;
\end{itemize}

\section{Equilateral Pentagons}
The equilateral pentagon in this study has the following attributes:
\begin{itemize} \setlength{\itemsep}{0mm}
\item Vertices are treated as point particles, each of unit mass;
\item Each edge is massless and fixed at unit length;
\item The angles between adjacent edges are allowed to change freely;
\end{itemize}
Note that the family of all equilateral pentagons includes degenerate pentagons
(e.g.\ an equilateral unit triangle with extra folded edges or a trapezium with one of the pentagonal angles taking $\pi$ ) and non-simple pentagons (e.g.\ the pentagram).

We denote the vertices of the pentagon by $z_i \in \C$ where $i \in \mathbb{Z}_5$, that is, the vertex indices are always modulo 5.
We will represent the elements of $\mathbb{Z}_5$ by $\{1,2,3,4,5\}$, starting with $1$.
The standard colour code that we use for the vertices are $\{z_1,z_2,z_3,z_4,z_5\}=\{$green, red, blue, black, yellow$\}$, this colour code will also be adopted for the relative angles $\psi_i$. The oriented edges of the pentagon are the vectors $z_{i+1} - z_i$.

A polygon is called {\em simple} if the edges do not intersect except at the vertices. 
The internal angle sum of a simple, $n$-sided polygon is $(n-2)\pi$, since its interior 
can be tiled by $n-2$ disjoint triangles.
This construction does not work when a polygon is self-intersecting. 
In fact, even the notion of an internal angle is not well defined in a self-intersecting polygon.
For this reason we adopt a convention for measuring the angles of the equilateral pentagon that gives the natural result of $\pm\frac{3\pi}{5}$ and $\pm\frac{\pi}{5}$, correspondingly, for all of the relative angles of the regular convex pentagons and regular pentagrams, where the sign denotes orientation.
We define the {\em relative angle} $\psi_i$ as the amount of rotation needed 
to turn the oriented edge $ z_{i+1} -  z_i  $ into the negative of the previous oriented edge $ z_{i} - z_{i-1}$, hence
$e^{i \psi_i} (z_{i+1} -  z_i ) = - ( z_{i} - z_{i-1})$ 
where $i \in \mathbb{Z}_5$.
As a result of this definition the relation between successive vertices $z_i \in \mathbb{C} $ is 
\begin{equation} \label{eq:zrecursion}
   z_{i+1} = z_i - e^{-i \psi_i} ( z_i - z_{i-1} ).
\end{equation}
For definiteness we use the principal argument $\Arg$ to define $\psi_i$ 
from the vertices $z_i$, such that
\begin{equation} \label{eq:psidef}
	\psi_i=-\Arg{\left(-\frac{z_{i+1}-z_i}{z_i-z_{i-1}}\right)}.
\end{equation}

The relative angles $\psi_i$ define the shape of the pentagon. 
The shape of the pentagon is invariant under rotations $z_i \to e^{i\theta} z_i$.
The absolute
angle $\theta$ is introduced to measure the orientation of the pentagon. 
We define $\theta$ to be the angle between the positive $x$-axis and the edge $z_2- z_1$ measured counter-clockwise from the $x$-axis, such that $z_2 = z_1 + e^{i \theta}$.

The shape of the pentagon is invariant under translations $z_i \to z_i + z$.
Initially we let the first vertex be arbitrarily located at $z$, but this $z$ will later be eliminated by fixing the centre of mass at the origin.
The vertices are:
\[
\begin{aligned}
\indent z_1&=z \\
\indent z_2&=z_1+e^{i \theta}\\
\indent z_3&=z_2-e^{i \theta} e^{-i\psi_2}\\
\indent z_4&=z_3+e^{i \theta} e^{-i(\psi_2+\psi_3)}\\
\indent z_5&=z_4-e^{i \theta} e^{-i(\psi_2+\psi_3+\psi_4)} \,.
\end{aligned}
\]
By a translation we can achieve $z_1 + z_2 + z_3 + z_4 + z_5 = 0$ so that the centre of mass is at the origin, and thus eliminating $z$ gives:
\begin{subequations}
\begin{equation}\label{eq:z1}z_1=\frac{1}{5} e^{i \theta } \left(-4+3 e^{-i \psi_2}-2 e^{-i (\psi_2+\psi_3)}+e^{-i (\psi_2+\psi_3+\psi_4)}\right)\end{equation}
\begin{equation}\label{eq:z2}z_2=\frac{1}{5} e^{i \theta } \left(1+3 e^{-i \psi_2}-2 e^{-i (\psi_2+\psi_3)}+e^{-i (\psi_2+\psi_3+\psi_4)}\right)\end{equation}
\begin{equation}\label{eq:z3}z_3=\frac{1}{5} e^{i \theta } \left(1-2 e^{-i \psi_2}-2 e^{-i (\psi_2+\psi_3)}+e^{-i (\psi_2+\psi_3+\psi_4)}\right)\end{equation}
\begin{equation}\label{eq:z4}z_4=\frac{1}{5} e^{i \theta } \left(1-2 e^{-i \psi_2}+3 e^{-i (\psi_2+\psi_3)}+e^{-i (\psi_2+\psi_3+\psi_4)}\right)\end{equation}
\begin{equation}\label{eq:z5}z_5=\frac{1}{5} e^{i \theta } \left(1-2 e^{-i \psi_2}+3 e^{-i (\psi_2+\psi_3)}-4 e^{-i (\psi_2+\psi_3+\psi_4)}\right).\end{equation}\label{eq:vertices}
\end{subequations}
In the following the symbol $z_i$ will refer to these formulas.

The ordered set of vertices $z_i$ defined by ~\eqref{eq:vertices} gives the equilateral pentagon modulo translations.
The ordered set of relative angles $\psi_i$ gives the equilateral pentagon modulo the special Euclidean group $SE(2)$
of orientation preserving rotations and translations. Thus the relative angles describe the labelled (and hence 
oriented) shape of the equilateral pentagon obtained by reducing the continuous symmetry $SE(2)$.
Later we will consider additional {\em discrete} symmetries to reduce further. 
They are the symmetry group $D_5$ of the equilateral (and equal masses) pentagon, and the 
reflection $\Z_2 = E(2)/ SE(2)$. The quotient by the full symmetry group $D_5 \times E(2)$ gives (unlabelled) shape up to congruence.

\subsection{Constraints}
When the vertices are placed at $z_i = i$ and consecutive edges are connected with joints at $z_2$, $z_3$ and $z_4$ 
we obtain the so called 4-segment open linkage system.
The shape space of this system is topologically a 3-dimensional torus, with angles $\psi_2$, $\psi_3$, and $\psi_4$.
The open linkage can be closed by requiring that $|z_1 - z_5| = 1$, introducing the fifth edge of length 1 connecting $z_5$ and $z_1$.
The closure constraint reduces the dimension of shape space to 2, and turns its topology into that of a rather 
complicated surface of genus 4 \cite{genus4a,genus4b}. 

The relative angles $\psi_5$ and $\psi_1$ do not enter the equations for $z_i$, see \eqref{eq:vertices},
and are completely determined by the other angles.
The equation for $\psi_5$ is obtained from \eqref{eq:zrecursion} for $i = 5$,
both sides are multiplied by $e^{i(\psi_2+\psi_3+\psi_4)}$, and then the complex logarithm is taken to obtain
\begin{equation}\label{eq:epsi5}
	\psi_5=-\Arg\left(1-e^{i \psi_4}+e^{i(\psi_3+\psi_4)}-e^{i (\psi_2+\psi_3+\psi_4)}\right).
\end{equation}
A similar calculation is used to show that
\begin{equation}\label{eq:epsi1}
	\psi_1=\Arg\left(1-e^{-i\psi_2}+e^{-i(\psi_2+\psi_3)}-e^{-i(\psi_2+\psi_3+\psi_4)}\right) \,.
\end{equation}
The above angles are all related by 
\begin{equation}
	\psi_1+\psi_2+\psi_3+\psi_4+\psi_5=(1+2k)\pi, \quad k \in \Z \,.
\end{equation} 
The sum of all five relative angles $\sum \psi_i$ takes values $-3\pi, -\pi, \pi, 3 \pi$, corresponding to $k = -2, -1, 0, 1$.
The extremal values are achieved in the region around the regular simple pentagon
where the pentagon remains convex, with either positive or negative orientation. 
Beyond these regions the pentagon may 
or may not be simple, but its angle sum remains constant as long as no additional 
stretched edge-configuration with $\psi_i = \pm \pi$ appears.
We could define the orientation of a pentagon by the sign of $\sum \psi_i$. 
For shapes with $\psi_i = \pm \pi$ this orientation is undefined as 
nearby shapes have either sign.

\begin{figure}[tp]
\centering
\subfigure[3 segments of the pentagon shown.]{
\includegraphics[width=7cm]{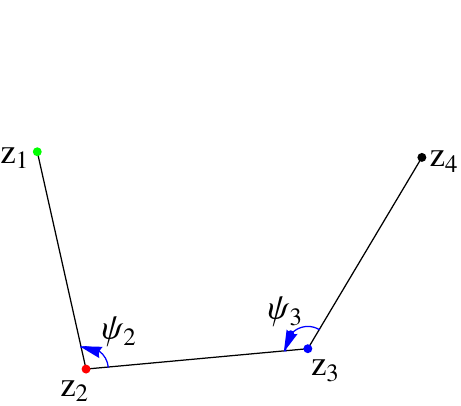}
\label{subfig:3seg}
}
\subfigure[With $\psi_2$ and $\psi_3$ specified, there are two possible pentagons in general.]{
\includegraphics[width=7cm]{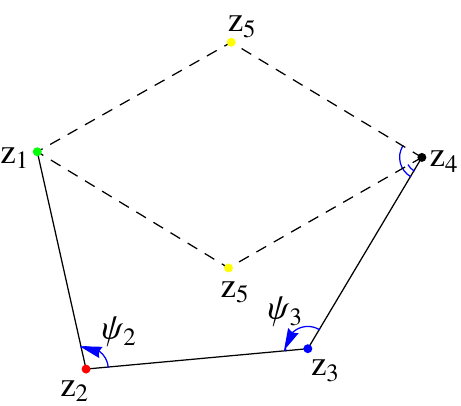}
\label{subfig:2pop}
}
\caption{Pentagon construction by specifying two successive relative angles.}
\label{fig:angles}
\end{figure}

The complicated topology of shape space of equilateral pentagons arises because the relative angles 
$\psi_2$, $\psi_3$, and $\psi_4$ are not independent.
Given two arbitrary angles $\psi_2$ and $\psi_3$, the relative orientation of three of the edges of the pentagon is fixed as shown in Figure \ref{fig:angles}\subref{subfig:3seg}. For the pentagon to be able to close with the two remaining segments, it is necessary that $0 \leq |z_4-z_1| \leq 2$. 
There are four possible cases when attempting to determine $\psi_4$ from $\psi_2$ and $\psi_3$:
\begin{itemize} \setlength{\itemsep}{0mm}
\item $|z_4-z_1|=0$: $\psi_4$ is undetermined, there are infinitely many solutions. 
\item $0<|z_4-z_1|<2$: the generic case where there are exactly two solutions for $\psi_4$ as illustrated in Figure \ref{fig:angles}\subref{subfig:2pop}.
\item $|z_4-z_1|=2$: the special case that there is a unique solution for $\psi_4$.
\item $|z_4-z_1|>2$: there is no solution since the pentagon cannot close.
\end{itemize}
The existence of the first case with infinitely many solutions makes the shape space 
more complicated than some gluing of the 2-torus $(\psi_2, \psi_3)$ of the  3-segment open linkage.
Starting with the closure condition
\begin{equation}\label{eq:congen}
	|z_5-z_1|^2-1=0,
\end{equation}
and inserting \eqref{eq:vertices} gives the relation between the relative angles 
$\psi_2$, $\psi_3$, and $\psi_4$ as
\begin{equation}
\label{eq:shapespacepsi}
\begin{split}
	3-2\cos{\psi_2}-2\cos{\psi_3}-2\cos{\psi_4}+2\cos{(\psi_2+\psi_3)} + \\
				+  2\cos{(\psi_3+\psi_4)}-2\cos{(\psi_2+\psi_3+\psi_4)}=0 \,.
\end{split}
\end{equation}
This equation defines the shape space as a 2-dimensional sub-manifold of the 3-dimensional torus,
the shape space of the 4-segment open linkage.
Equation \eqref{eq:congen} can be rewritten as $X \zeta +  \zeta \bar\zeta +X^{-1}  \bar \zeta = 0 $ where 
$X = e^{i \psi_4}$, $\zeta = -1 + e^{i \psi_3} - e^{i(\psi_2 + \psi_3)}$, such that 
$\zeta \bar\zeta = |z_4 - z_1|^2 = 3-2\cos{\psi_2}-2\cos{\psi_3}+2\cos{(\psi_2+\psi_3)}$.
Assuming $|\zeta| \not = 0$ , and using the polar form of $\zeta$ gives
\begin{equation}
	\psi_4 = \pm \arccos{\left(-\frac{|\zeta|}{2}\right)} -  \arg{\zeta}  \mbox{\hspace{1cm}s.t. } \psi_4 \in (-\pi,\pi]\,. 
\end{equation}
When $|\zeta| =  |z_4 - z_1| = 0$ the value of $\psi_4$ is undetermined. 
This occurs only when $\psi_2  = \psi_3 = \pm \pi/3$.
Thus, solving of \eqref{eq:congen} reflects exactly the four cases discussed above.

\subsection{Shape Space}
Although the relative angles $\psi_i$ are easy to visualise and interpret geometrically, the algebraic equations can be simplified with the following affine transformation
\begin{equation}\label{eq:psi2alpha}
\left[{\begin{array}{c}\alpha_1\\\alpha_2\\\alpha_3\end{array} }\right] = 
\left[{\begin{array}{rrr}
-\frac{1}{2} & -1 & -\frac{1}{2}\\
-\frac{1}{2} & 0 & -\frac{1}{2}\\
-\frac{1}{2} & 0 & \frac{1}{2}\\\end{array} }\right]
\left[{\begin{array}{c}\psi_2\\\psi_3\\\psi_4\end{array} }\right]  + 
\left[{\begin{array}{c}0\\\pi\\0\end{array} }\right]
= A \vec{\psi} + \textbf{b}.
\end{equation}

The determinant of $A$ is $-\frac{1}{2}$, the negative sign means that the orientation is reversed, and the factor $\frac{1}{2}$ means that the area is halved when we do the transformation. That means if we use the same natural domain $(-\pi, \pi]$ for both the $\psi$ and $\alpha$-coordinates, then we have a double covering in the $\alpha$-coordinates. That is, 
every possible pentagon in the $\psi$-coordinates occurs exactly twice in the $\alpha$-coordinates. 
Since the angles $\psi_i$ are defined modulo $2\pi$ this induces an equivalence relation for the angles $\alpha_i$, which 
is given by
\begin{equation} \label{alphaequiv}
   (\alpha_1, \alpha_2, \alpha_3) \equiv 
   (\alpha_1, \alpha_2, \alpha_3) 
   - (  (i + 2 j + k) \pi,  (i + k) \pi, ( i - k) \pi)
    \qquad \text{for}\quad  i, j, k \in \Z \,.
\end{equation}
Using this equivalence relation the double covering can be removed by restricting the fundamental domain to our choice of $\alpha_3 \in [0,\pi)$.

The transformation is invertible, and its inverse is
\begin{equation}\left[{\begin{array}{c}\psi_2\\\psi_3\\\psi_4 \end{array} }\right] = \left[{\begin{array}{rrr}
0 & -1 & -1 \\
-1 & 1 & 0 \\
0 & -1 & 1 \\
\end{array} }\right]
\left[{\begin{array}{c}\alpha_1\\\alpha_2  -\pi \\\alpha_3\end{array} }\right]  \bmod 2 \pi.
\end{equation}

\begin{figure}[tb]
\centering
\includegraphics[width=9cm]{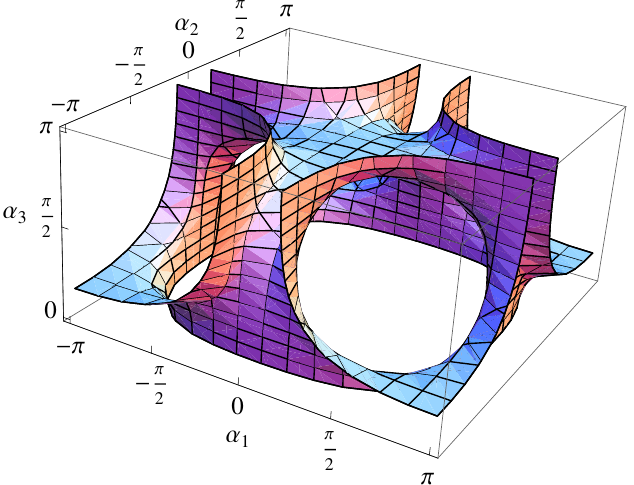}
\caption{Shape space in $\alpha$-coordinate system. 
In $\alpha_1$ and $\alpha_2$ there are periodic boundary conditions.
The gluing from $\alpha_3  =0$ to $\alpha_3 = \pi$ is
done after a shift by $\pi$ in $\alpha_1$ and $\alpha_2$.
This is a surface of genus 4.}\label{fig:shapespace}
\end{figure}

For convenience, both coordinate systems will be used. The geometrical interpretations will be done in the $\psi$-coordinate system while the algebraic calculations will be done in the $\alpha$-coordinate system.
In $\alpha$-coordinates the equation for shape space \eqref{eq:shapespacepsi} simplifies to the symmetric form
\begin{equation}\label{eq:conalpha}
	C(\alpha_1, \alpha_2, \alpha_3) := 3+4\cos{\alpha_1}\cos{\alpha_2}+4\cos{\alpha_1}\cos{\alpha_3}+4\cos{\alpha_2}\cos{\alpha_3}=0.
\end{equation}
The corresponding surface $C(\alpha_1, \alpha_2, \alpha_3) = 0$ is shown in Figure~\ref{fig:shapespace}. 
It is reminiscent of the I-WP surface from minimal surface theory \cite{Hyde97}.
Our defining equation is like a low order Fourier approximation of the I-WP surface.

From equation \eqref{eq:conalpha}, $\alpha_3$ can be uniquely expressed in terms of $\alpha_1$ and $\alpha_2$ as
\begin{equation}\label{eq:alpha3equals}
	\alpha_3=\arccos{\left(\frac{-3-4\cos{\alpha_1}\cos{\alpha_2}}{4(\cos{\alpha_1}+\cos{\alpha_2})}\right)}.
\end{equation}
A unique solution is obtained for $\alpha_3$ as we are strictly in the domain $[0,\pi)$.
Accordingly $\alpha_1$ and $\alpha_2$ are good local coordinates almost everywhere on shape space. The exceptions are at the boundaries where $\alpha_3 = 0, \pi$, and, 
in particular, 
when $(\alpha_1, \alpha_2)$ takes the values $(\pm \pi/6, \pm 5 \pi/6)$ and $(\pm 5 \pi/6, \pm \pi/6)$.
Over these points there is a whole line in the surface corresponding to
infinitely many solutions of $\alpha_3$. This corresponds to the case when 
the fraction in \eqref{eq:alpha3equals} is of the form $0/0$.
All other points for which $\cos\alpha_1 + \cos\alpha_2=0$ are not on shape space at all, so we do not need to worry about them.
We will see that after discrete symmetry reduction it is not necessary to consider the singular lines anymore (where the exceptions for $\alpha_1$ and $\alpha_2$ occur).

Perhaps the simplest way to understand why the genus is 4 is to compute the Euler-characteristic
from the picture taking the gluing into account.
Alternatively we can compute it using Morse theory from the critical points of a smooth function defined on the 
surface, we will follow this approach with the moment of inertia later on.

Notice that our definition of shape space is the set of all equilateral pentagons with labelled (or distinguishable) vertices.
We use the colours green, red, blue, black, yellow to designate the vertices $z_1, z_2, z_3, z_4, z_5$, respectively.
In the standard notion of congruence of polygons the vertices (and sides) are unlabelled and hence indistinguishable.
Considering the polygon as a mechanical system (with potentially different masses at the vertices and different 
moments of inertia for the edges) we obtain a description of the space of labelled equilateral pentagons.
As it turns out, all quantities we are interested in here are given by functions on the unlabelled shape space, 
because the masses at the corners and the side-length are all equal. 
The reduction from labelled to unlabelled shapes is a discrete symmetry reduction, which we are going to study next.

\section{Discrete Symmetries}
Consider an arbitrary equilateral pentagon with 5 labelled vertices, e.g.\ distinguished by colours.
The action of the group $D_5$ on the labelled shape space is generated by
{\em vertex rotations} and {\em vertex reflections}.
By a vertex rotation we mean a cyclic permutation of the vertices, that is
$R(z_1, z_2, z_3, z_4, z_5) = (z_2, z_3, z_4, z_5, z_1)$.
By a vertex reflection we mean a permutation that reverses the order of vertices and fixes a chosen vertex; 
choosing to fix $z_3$, we get
$V(z_1, z_2, z_3, z_4, z_5) = (z_5, z_4, z_3, z_2, z_1)$.
Both operations leave the set of vertices $z_i$  in the plane fixed; they merely permute the vertices in a way that preserves neighbours. 
Thus the unlabelled shape of the pentagon is fixed as well, but the orientation induced by the labelling is reversed under $V$.
For a generic shape both operations do change the labelled shape, that is there is no element in $E(2)$
that undoes the action. 
Special labelled shapes are fixed under subgroups of $D_5$, for example the labelled regular pentagons are
invariant under vertex rotations $R$, since up to a geometric rotation $\in SE(2)$ it is the same labelled shape as before vertex rotation.
Symmetry reduction of the labelled shape space allows us to define a fundamental region in which each 
unlabelled shape is represented exactly once.
Since the action of $D_5$ is not free, the unlabelled shape space is not a smooth manifold but just an orbifold,
with singularities at the shapes that have higher symmetry, and hence non-trivial isotropy.

If we write the relative angles as a vector $\vec{\psi}=(\psi_1, \psi_2, \psi_3, \psi_4, \psi_5)^t$,
then the symmetry operations can be represented by multiplication of $\vec{\psi}$ by matrices $R$, for 
vertex rotation, and $V$, for vertex reflection (fixing vertex 3) where
\[
R=\left( \begin{array}{ccccc}
0 & 1 & 0 & 0 & 0\\
0 & 0 & 1 & 0 & 0\\
0 & 0 & 0 & 1 & 0\\
0 & 0 & 0 & 0 & 1\\
1 & 0 & 0 & 0 & 0\end{array} \right), 
\quad
V=\left( \begin{array}{ccccc}
0 & 0 & 0 & 0 & -1\\
0 & 0 & 0 & -1 & 0\\
0 & 0 & -1 & 0 & 0\\
0 & -1 & 0 & 0 & 0\\
-1 & 0 & 0 & 0 & 0\end{array} \right) \,.
\]
The additional minus sign in $V$ follows from the definition of $\psi_i$ in \eqref{eq:psidef}.
The group $D_5$ is generated by $R$ and $V$ with presentation 
$\langle R,V \,|\, R^5=V^2=id,VR=R^{-1}V \rangle $.
The group $D_5$ acting on the relative angles $\vec{\psi}$ leaves the unlabelled (and un-oriented) shape invariant.
In general it changes the labelled shape; normally there are 10 (the order of $D_5$) different labelled shapes corresponding to the same unlabelled un-oriented shape. 

%

There is another discrete symmetry because by using the relative angles $\psi_i$ to describe the shape 
we have reduced by $SE(2)$, but not by $E(2)$. Hence reflections about a line through the origin
give another discrete symmetry $\Z_2 = E(2) / SE(2)$.
We call it the mirror reflection symmetry $M$. The action on the vertices is 
$M(z_1, z_2, z_3, z_4, z_5) = (\bar z_1, \bar z_2, \bar z_3, \bar z_4, \bar z_5)$, 
where the overbar denotes complex conjugation.
As a matrix acting on the space of angles $\vec{\psi}$ we simply have $M = -id$.
Unlike $R$ and $V$ the operation $M$, in general, changes the set of vertices $z_i$,
even modulo $SE(2)$.

Labelled polygons have an orientation induced by the labelling, 
while for unlabelled polygons an orientation may be kept track of by orienting the edges with an arrow. 
Considering oriented pentagons, both $M$ and $V$ reverse the orientation.
Simply forgetting the labels of the vertices gives an unlabelled un-oriented shape.
This corresponds to reduction by the full group $D_5 \times E(2)$, which gives a shape in the 
classical sense that two shapes are the same if they are congruent (and have the same size).
We will see that unlabelled {\em oriented} shapes up to congruence are obtained from reduction 
by a slightly different group $D_5^+ \times SE(2)$.

Combining all three discrete symmetries gives  $D_5 \times \Z_2$ which is isomorphic to 
$D_{10}$ with presentation $\langle MR,V \,|\, (MR)^{10}=V^2=id,V(MR)=(MR)^{-1}V \rangle$.
A subgroup of  $D_{10}$, different (but isomorphic) to $D_5$ generated by $R$ and $V$, 
is obtained from the generators $R$ and $MV$:
$D_5^+ = \langle R, MV \,|\, R^5=(MV)^2=id,(MV)R=R^{-1}(MV) \rangle $.
The superscript $+$ indicates that the action of this group preserves orientation as both $R$ and $MV$ preserve orientation.
The group $D_5^+$ is used for the symmetry reduction in the next chapter, and the resulting reduced shapes
are unlabelled oriented shapes. In the following we will drop the qualification unlabelled and simply 
talk about oriented shapes.

\subsection{Discrete Symmetry Reduction}
We now construct a fundamental region of labelled shape space
such that the whole surface is obtained as the $D_{5}^+$ orbit of this fundamental region.
This fundamental region will contain every oriented shape exactly once.
If necessary, a final reduction by $V$ to remove the double covering of un-oriented shapes from this region may be performed to halve the fundamental region.

The important objects in symmetry reduction are  isotropy subgroups.
The isotropy group of a point $x$ in shape space is defined as $G_x = \{ g \in G \, : \, g x = x \}$.
A shape has high symmetry if it has a large isotropy group.
In simple cases it is enough to consider the length of the orbit of $x$, namely  $\# \{ g x \, : \, g \in G \}$
to distinguish different isotropy types. 
We now discuss the isotropy groups that occur for the action of $D_5^+$ on labelled shapes.

The shapes with the highest symmetry are regular pentagons, which are the pentagrams $\vec{\psi} = \pm(\frac{\pi}{5},\frac{\pi}{5},\frac{\pi}{5},\frac{\pi}{5},\frac{\pi}{5})^t$ 
and the regular convex pentagons $\vec{\psi} = \pm(\frac{3\pi}{5},\frac{3\pi}{5},\frac{3\pi}{5},\frac{3\pi}{5},\frac{3\pi}{5})^t$.
These special shapes are fixed under the whole group $D_5^+$, that is their isotropy group is $D_5^+$, and their
orbit length is 1.

The next group of symmetric labelled shapes are reflection symmetric, with respect to some axis through the centre of mass, which is fixed at the origin.
After relabelling the vertices by some power of $R$ the relative angles of a symmetric shape are $(\psi_1, \psi_2, \psi_3, \psi_2, \psi_1)$,
which is clearly fixed under $MV$.
Hence the isotropy group of these shapes is $\Z_2$ generated by $MV$ (or $R^k MV R^{-k}$), and their orbit length is 5.
Examples of reflection symmetric shapes are shown in Figure~\ref{fig:time} at $t = 0, \tau/10, \tau$.

All other labelled shapes have a trivial $D_5^+$ isotropy group, so their orbit under $D_5^+$ is length 10, and they do not posses any symmetry.

The fundamental region of the action of $D_5^+$  is constructed in the $\alpha$-space
by using the reflection symmetric shapes as a boundary.
The action of $V$ and $M$ on $(\alpha_1, \alpha_2, \alpha_3)$ is simple: 
$V ( \alpha_1, \alpha_2, \alpha_3) =  (-\alpha_1, -\alpha_2, \alpha_3)$
and 
$M ( \alpha_1, \alpha_2, \alpha_3) =  (-\alpha_1, -\alpha_2, -\alpha_3) 
\equiv ( \pm \pi - \alpha_1, \pm \pi - \alpha_2, \pi - \alpha_3)$
where the $\pm$ sign is negative if $\alpha_i$ is negative, and positive otherwise, 
then by using the equivalence relation \eqref{alphaequiv} we get
$MV ( \alpha_1, \alpha_2, \alpha_3) = (\alpha_1, \alpha_2, -\alpha_3) 
\equiv (\pm \pi + \alpha_1, \pm \pi + \alpha_2, \pi - \alpha_3)$.
Thus labelled shapes with $\alpha_3 = 0$ are fixed under $MV$.
The operation $R$ in explicit terms is somewhat more complicated, but we do not 
require that formula.

Note that $\alpha_3 = 0$ implies $\psi_2=\psi_4$ and $\psi_1 = \psi_5$, which is an expression 
of the reflection symmetry with respect to vertex 3.
Each point on the curve $\alpha_3 = 0$ is fixed under $MV$, while $M$ and $V$ leave the curve
invariant as a whole. In fact the action of $M$ and $V$ on the curve $\alpha_3 = 0$ are the same, 
simply $(\alpha_1, \alpha_2) \to (-\alpha_1, -\alpha_2)$.
We denote the curve $\alpha_3 =0$ and its images under $D_5^+$ as  {\em symmetry curves}.
These curves are fixed sets of involutive elements in $D_5^+$ of the form $R^k MV R^{-k}$, 
where $k=0$ corresponds to the symmetric shapes for which vertex $z_3$ is on the symmetry line,
since $MV$ fixes $\psi_3$.

\begin{figure}[tbp]
\centering
\includegraphics[width=7cm]{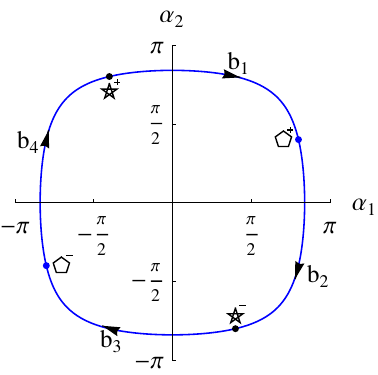}
\caption{The basic symmetry curve $\alpha_3 = 0$ that is used to generate the division of the labelled shape space. One quarter of it, denoted by $b_1$, can be transformed into the curve $b_3$ by the symmetry operation $V$ or $M$ and vice versa; similarly with the curves $b_2$ and $b_4$. Notice how all the $b_i$ curves are $\frac{\pi}{2}$ rotations of each other. 
The superscript $\pm$ signs on the pentagons indicate orientation.}\label{fig:Bcurves}
\end{figure}

The orbit of reflection symmetric shapes under $D_5^+$ has length 5, and hence there are 5 symmetry curves
obtained by letting $R$ the generator of $D_5^+/\Z_2$ act on the basic symmetry curve $\alpha_3 = 0$.
These 5 curves divide the surface of labelled shapes into 10 simply connected regions, 
since 10 is the length of the orbit under $D_5^+$ of a generic (that is non-symmetric) labelled shape.

Explicit formulas for the basic symmetry curve $\alpha_3 = 0$ are obtained from $C(\alpha_1, \alpha_2, 0) = 0$.
Figure~\ref{fig:Bcurves} shows a quarter of the basic symmetry curve parametrized as
\begin{equation}\label{eq:b1}
b_1(t)=\left(t,\arccos{\left[-\frac{1}{8}(3+4\cos{t})\sec^2{\left(\frac{t}{2}\right)}\right]},0\right)_\alpha \mbox{ where } t=\left[-\frac{2\pi}{5}, \frac{4\pi}{5}\right)
\end{equation}
and the subscript $\alpha$ is used to denote the $\alpha$-coordinate system.

\begin{figure}[tb]
\centering
\includegraphics[width=10cm]{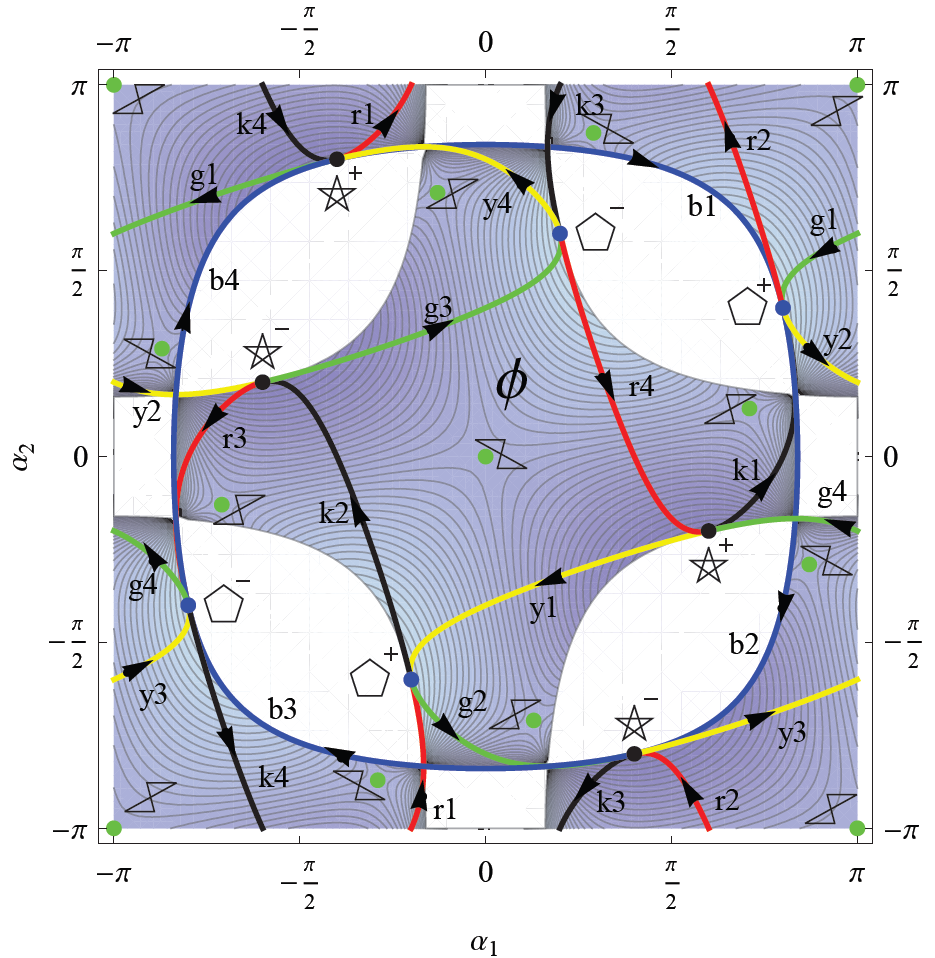}
\label{subfig:symcur2d}
\caption{Projection of the shape space onto the plane $(\alpha_1, \alpha_2)$ with contours of constant moment of inertia. 
Taking the rotation symmetry $R^k$ for each of the blue curves $b_i$ from Figure \ref{fig:Bcurves} produces the red, green, yellow, black curves denoted by $r_i$, $g_i$, $y_i$, $b_i$, respectively, 
for $i = 1,2,3,4$. These are the fixed sets of $R^k MV R^{-k}$ for $k = 0, \dots, 4$ with vertex $z_{3-k}$ fixed when defining the vertex reflection.
Of the 10 topologically equivalent regions we choose the piece bounded 
by the curves $r_4, g_3, k_2$ and $y_1$ as our fundamental region and denote it by $\phi$.
}\label{fig:symcur}
\end{figure}

When $t=-\frac{2\pi}{5}$, the relative angles are $(-\frac{2\pi}{5},\frac{4\pi}{5},0)_\alpha=(\frac{\pi}{5},\frac{\pi}{5},\frac{\pi}{5},\frac{\pi}{5},\frac{\pi}{5})^t$, this gives the positively oriented pentagram.\footnote{$(\alpha_1,\alpha_2,\alpha_3)_\alpha=(\psi_1,\psi_2,\psi_3,\psi_4,\psi_5)^t$ shows equivalent relative angles in different coordinate systems.} 
When $t=\frac{4\pi}{5}$, the relative angles are $(\frac{4\pi}{5},\frac{2\pi}{5},0)_\alpha=(\frac{3\pi}{5},\frac{3\pi}{5},\frac{3\pi}{5},\frac{3\pi}{5},\frac{3\pi}{5})^t$, this corresponds to the positively oriented regular convex pentagon.
Applying $M$ or $V$ to $b_1$ gives $b_3$. In the plane $\alpha_3 = 0$ this amounts is a rotation by $\pi$. 
If instead a rotation by $\frac{\pi}{2}$ is performed the curves $b_2$ and $b_4$ are obtained.
This is a result of the fact that the equation $C(\alpha_1, \alpha_2, \alpha_3) = 0$ in \eqref{eq:conalpha} is even in $\alpha_i$.

The images of the basic symmetry curve of shapes for which vertex $z_3$ is on the line of reflection symmetry in the plane 
are shown in Figure \ref{fig:symcur}. The basic symmetry curves are the fixed set of the involution $MV$.
The $k^{th}$ image of this curve under $R$ is the fixed set of $R^k MV R^{-k}$, with vertex $z_{3-k}$ on the line of 
reflection symmetry in the plane.
The colour scheme for vertices and symmetry curves was chosen, such that shapes on the symmetry curve with colour $c$,
have the vertex with that colour on the line of reflection symmetry in the plane.

\begin{figure}[tb]
\centering
\includegraphics[width=11.5cm]{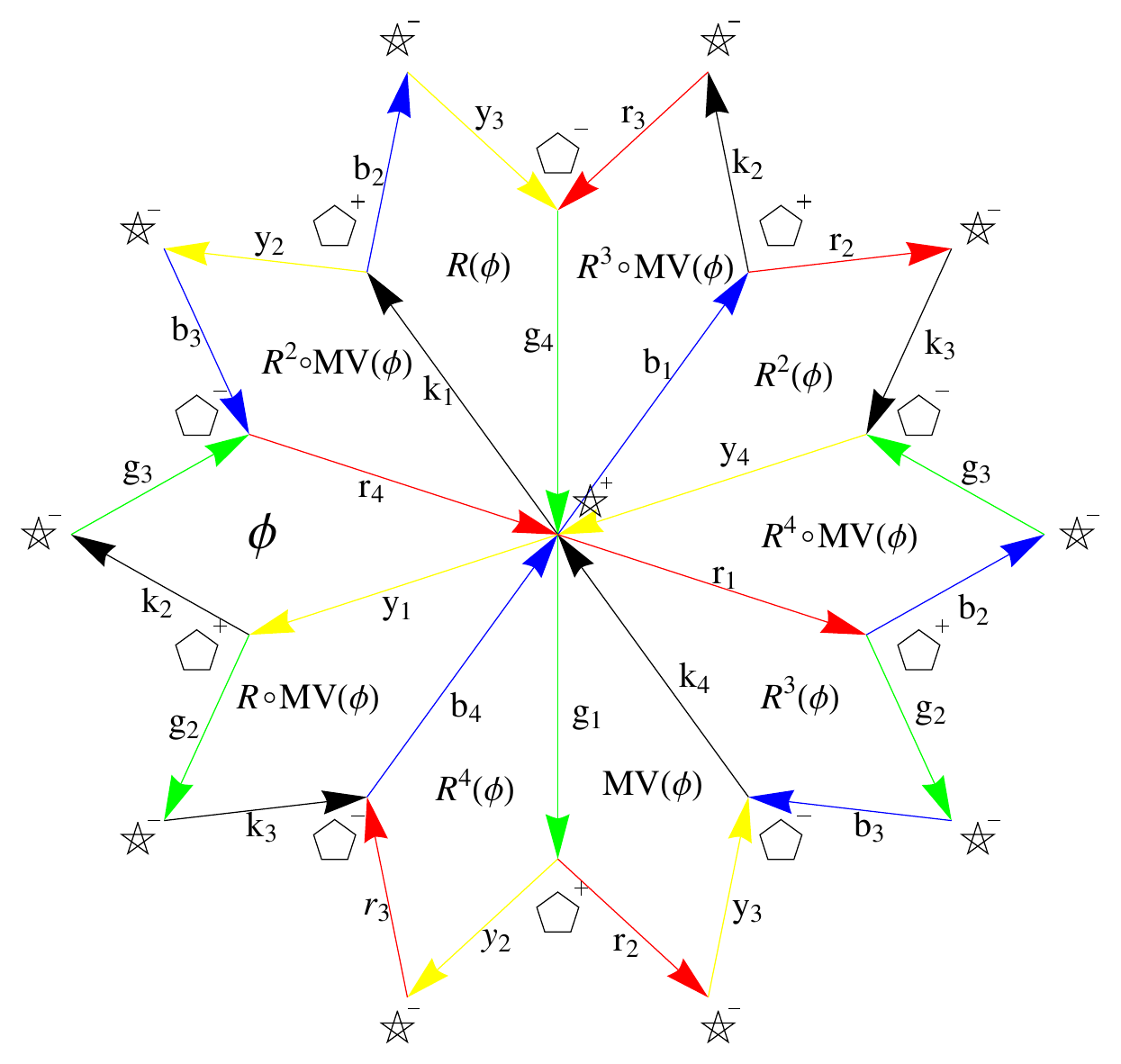}
\caption{Schematic picture of how to glue the labelled shape space from copies of the fundamental region $\phi$. The labelled outer edges with their specified orientation are to be glued together to form the shape space of genus 4. The symmetry operations required to map $\phi$ to any other region is identified in this diagram.
}\label{fig:map}
\end{figure}

We chose a fundamental region $\phi$ as shown in Figure \ref{fig:symcur}, which is bounded by 
$R(b_4)$, $R^3(b_1)$, $R^4(b_2)$, and $R^2(b_3)$.
The fundamental region $\phi$ has the useful property that it can be parameterised by $\alpha_1$ and $\alpha_2$.
Starting with $\phi$, the nine remaining regions are obtained by applying the symmetries $R$ and $MV$ to $\phi$.
When considering an arbitrary point in shape space, it may always be moved to the fundamental region $\phi$ by applying some symmetry operation
from $D_5^+$.
The topology of how the 10 pieces are glued together to give all of the shape space is shown in Figure~\ref{fig:map}.
This figure is somewhat reminiscent of a tiling of the hyperbolic plain obtained from the triangle group of type $(2,5,5)$, 
but since the tiles are quadrilaterals in our case the analogy is rather incomplete.
Now all the analysis required can be done in $\phi$. The remaining parts of the surface are covered by the $D_5^+$ orbit of $\phi$.

In Figure \ref{fig:symcur} (and later figures) two copies of each orientation of the regular pentagon and pentagram are shown. 
Note that these are equivalent under \eqref{alphaequiv}, so they represent the same point of the labelled shape space.
The reason for displaying both is to emphasise the peculiarity of our fundamental region, which occurs because the corner points of the fundamental region $\phi$ are actually outside of the fundamental domain of the $\alpha's$ since $\alpha_3 \in [0, \pi)$.


Finally, we discuss the reduction by the remaining symmetry $V$, which in general reverses the orientation of a given oriented shape.
Our $D_5^+$ fundamental region $\phi$ contains both of these, since they are different as oriented shapes.
The quotient by $V$ identifies these oriented shapes with different orientations as the same, 
and the result is the set of un-oriented shapes, or just shapes in standard terminology.
A natural way to define a fundamental region that contains each (un-oriented) shape once 
(up to congruence $E(2)$) is to cut $\phi$ into halves along $\alpha_1 = \alpha_2$. 
Along this line the total angle sum $\sum \psi_i$ jumps from $\pi$ to $-\pi$.
We denote $\phi^+$ as the half of $\phi$ that contains each shape with a positive angle sum.
Thus $\phi^+$ contains each (un-oriented) shape exactly once.
Half of the additional side $\alpha_1 = \alpha_2$ of $\phi^+$ is open while the other half, including the midpoint, is closed. Every point in the interior of $\phi^+$ is generic and has a $D_{10}$ orbit of length 20. The origin has orbit length 10, while the sides have length 10 and the corners have length 2.
Hence there is one new symmetric shape, which is invariant under $V$ (but not under $M$), corresponding to the 
origin $\alpha_1 = \alpha_2 = 0$ of $\phi$. This peculiar shape is shown in Figure~\ref{subfig:moisada} and is invariant under orientation reversal.

\section{Moment of Inertia and Angular Momentum}
The moment of inertia of the equilateral pentagon with point masses $m_i = 1$ at the vertices
with respect to its centre of mass is given by
\begin{equation}\label{eq:moiform}
I=\displaystyle\sum_{i=1}^5{|z_i|^2},
\end{equation}
which is obviously invariant under the symmetry group $D_{10}$.
In the $\alpha$-coordinate system equation \eqref{eq:moiform} becomes
\begin{equation}\label{eq:moialpha}
I=4+2\cos{\alpha_1}\cos{\alpha_2}+\frac{8}{5}\cos{\alpha_1}\cos{\alpha_3}+\frac{12}{5}\cos{\alpha_2}\cos{\alpha_3}+\frac{6}{5}\sin{\alpha_1}\sin{\alpha_2}\,.
\end{equation}
The contours of constant $I$ are shown in Figure~\ref{fig:symcur}.
Notice how \eqref{eq:moialpha} is even in $\alpha_3$; this implies $\alpha_3$ can be eliminated using \eqref{eq:alpha3equals} to give an expression that is a rational function of trigonometric functions.
The denominator of this function of two variables vanishes at the 8 points given by $(\pm \pi/6, \pm 5 \pi/6)_\alpha$ and $(\pm 5 \pi/6, \pm \pi/6)_\alpha$ where taking all possible combinations of the $\pm$ sign; this has already been
discussed in relation to \eqref{eq:alpha3equals}. 
These points are all outside the fundamental region so they do not cause a problem.
The moment of inertia is invariant under the full symmetry group $D_{10}$.
To find all critical points of the moment of inertia we employ symmetry reduction, and hence considered only the critical points of the moment of inertia within the fundamental region $\phi^+$. 
The action of $D_{10}$ then generates all critical points on the entire surface. 
On $\phi$, as shown in Figure \ref{fig:symcur}, there are five critical points of the moment of inertia. Four of them are on the boundary, and are in fact on the four corners and one is in the interior. 
The four critical points on the boundary of the fundamental region also happen to be located on the boundary $\alpha_3 = \pi$ where 
$(\alpha_1, \alpha_2)$ are not local coordinates. Thus the criticality is established using \eqref{eq:moialpha} and Lagrange multipliers 
incorporating the constraint \eqref{eq:conalpha}. 

\begin{figure}[tbp]
\centering
\subfigure[2 minima;  \newline $I = \frac{1}{2}(5-\sqrt{5})$;\newline isotropy: $\langle R, MV \rangle$]{
\includegraphics[width=3cm]{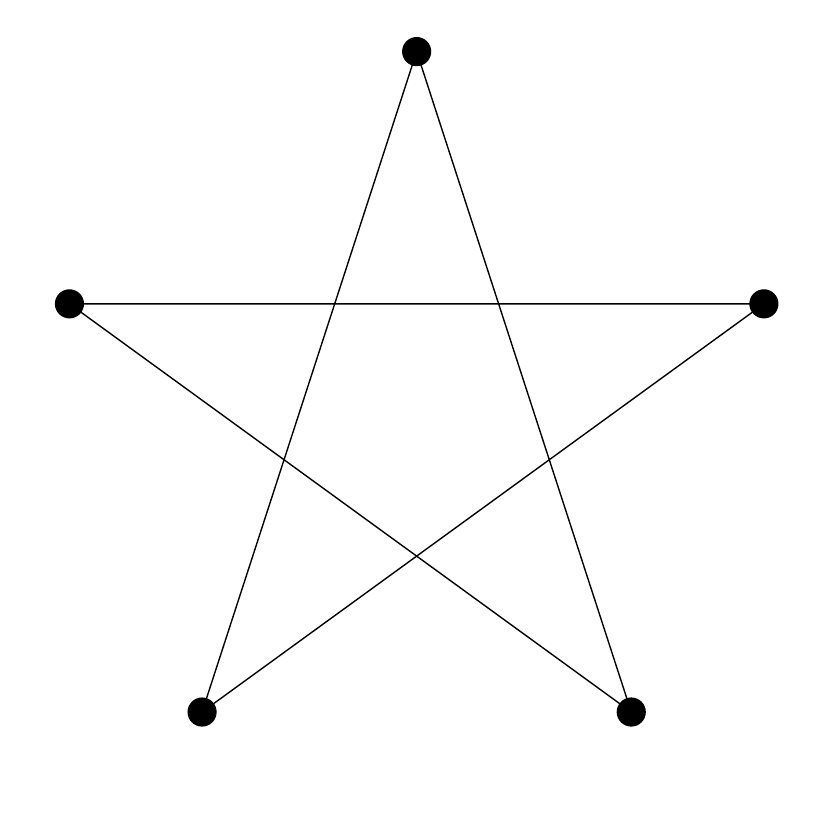}
\label{subfig:moimin}
}
\subfigure[2 maxima;  \newline $I = \frac{1}{2}(5+\sqrt{5})$;\newline isotropy: $\langle R, MV\rangle$]{
\includegraphics[width=3cm]{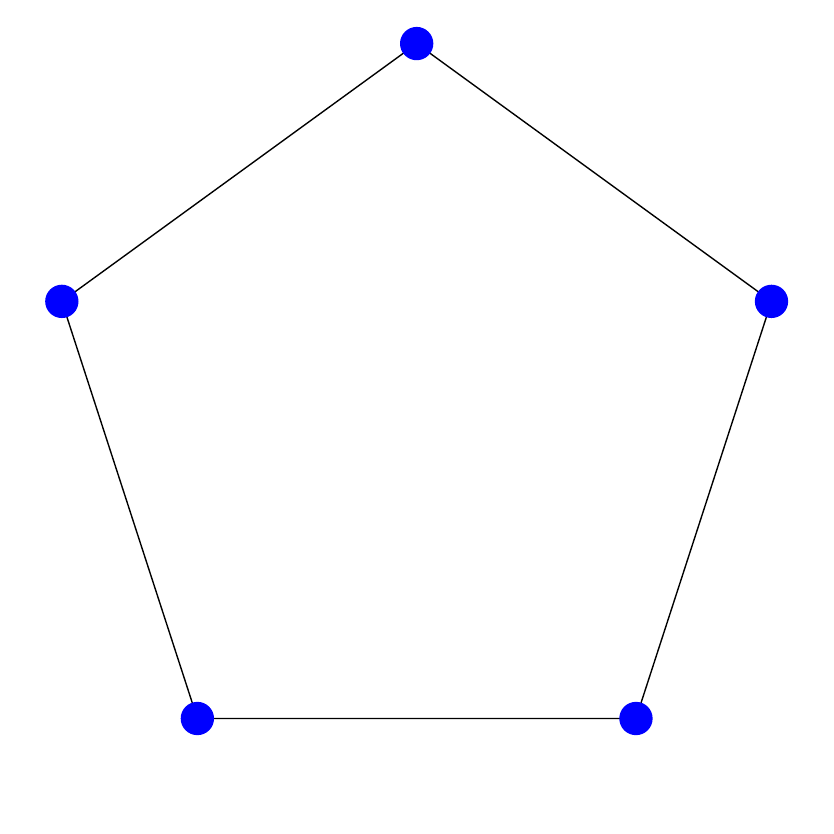}
\label{subfig:moimax}
}
\subfigure[5 saddles; 
\newline$I = \frac{5}{2}$;\newline isotropy: $\langle V \rangle$]{
\includegraphics[width=3cm]{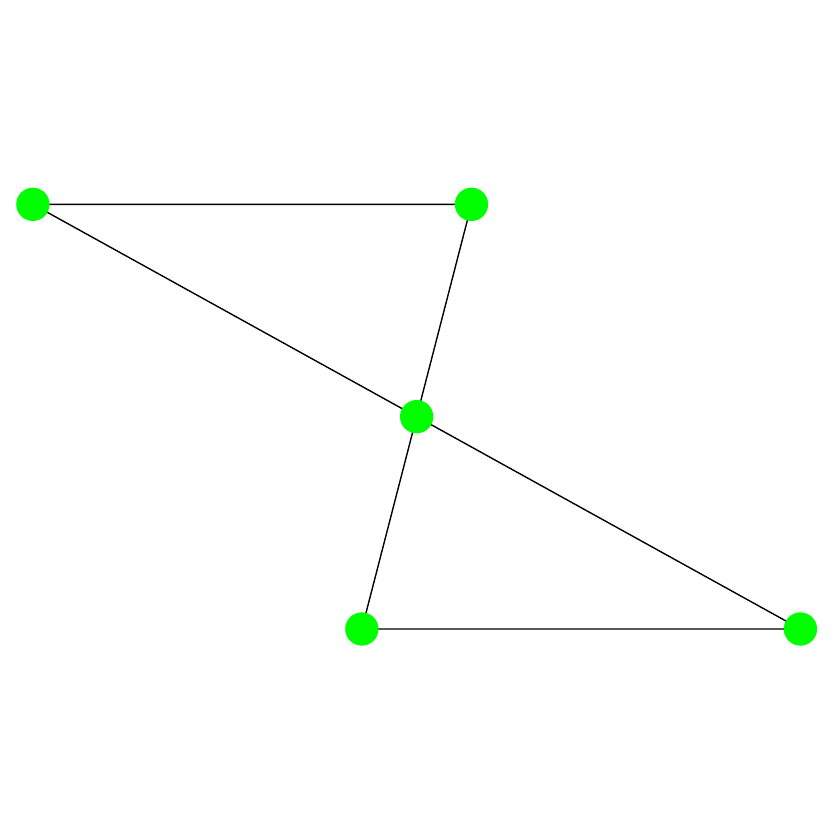}
\label{subfig:moisada}
}
\subfigure[5 saddles;\newline $I = \frac{5}{2}$;\newline isotropy: $\langle V \rangle$]{
\includegraphics[width=3cm]{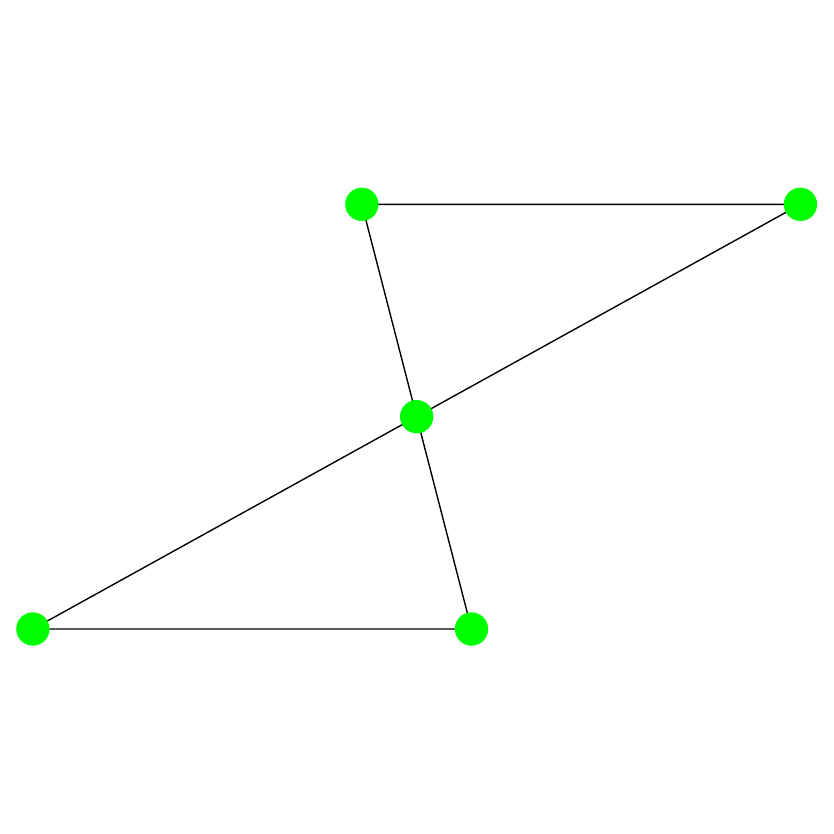}
\label{subfig:moisadb}
}
\caption{
Applying the symmetry operations $R, M$ and $V$ to the critical shapes shows that there are in total 14 critical labelled shapes. The types are: 2 minima [\ref{subfig:moimin}], 2 maxima [\ref{subfig:moimax}], 10 saddles (5 of the form [\ref{subfig:moisada}] and their mirror images under $M$ [\ref{subfig:moisadb}]). 
}\label{fig:moilist}
\end{figure}

The pentagrams in Figure \ref{subfig:moimin} have global minimal moments of inertia with value $\frac{1}{2}(5-\sqrt{5})$ and are located at $\pm(\frac{3\pi}{5},-\frac{\pi}{5},\pi)_\alpha=\pm(\frac{\pi}{5},\frac{\pi}{5},\frac{\pi}{5},\frac{\pi}{5},\frac{\pi}{5})^t$. 
The regular convex pentagons in Figure \ref{subfig:moimax} have global maximal moments of inertia equal to $\frac{1}{2}(5+\sqrt{5})$ and are located at $\pm(\frac{\pi}{5},\frac{3\pi}{5},\pi)_\alpha=\pm(\frac{3\pi}{5},\frac{3\pi}{5},\frac{3\pi}{5},\frac{3\pi}{5},\frac{3\pi}{5})^t$.
%
%
The saddle critical shapes shown in Figure \ref{subfig:moisada} and \ref{subfig:moisadb} have moments of inertia equal to $\frac{5}{2}$; and the saddle found in $\phi$ occurs at
$(0,0,\pi - \kappa)_\alpha = 
  \left(\frac{1}{2}( \pi - \kappa) ,\kappa,\pi,-\kappa,-\frac{1}{2}(\pi - \kappa)\right)^t$
where $\kappa = \arccos( \frac{7}{8} )$. Applying the symmetry operations to this saddle produces a total of 10 saddles on the entire surface.
Recalling the $D_{10}$ orbit lengths of the symmetric shapes there are 2 minima, 10 saddle points, and 2 maxima, 
so that the Euler characteristic of the shape space is 
$\chi = 2 - 10 + 2 = -6$ proving that the genus is 4, as claimed earlier.

For a system of particles with coordinates $(x_i, y_i)$ in the plane
the angular momentum with respect to the origin (which is the centre of mass in our case) is a scalar given by
\begin{equation}
L=\displaystyle\sum_{i=1}{m_i(x_i\dot{y_i}-\dot{x_i}y_i)}
\label{eq:angularmomentum2d}
\end{equation}
where the dot denotes the time-derivative, so that $\dot x_i$ is the $x$-component of the velocity vector of vertex $z_i$.
In the case of the equilateral pentagon where it is placed on the complex plane with each vertex having unit mass, equation \eqref{eq:angularmomentum2d} simplifies to 
\begin{eqnarray}
L 
&=& \displaystyle\sum_{i=1}^5{\Im{(\bar{z_i}\dot{z_i})}} \mbox{\indent\indent where } z_i=x_i+iy_i \, .\label{eq:angmomform}
\end{eqnarray}
Computing the angular momentum by transforming the vertices found from equation \eqref{eq:vertices} into the $\alpha$-coordinates and using equation \eqref{eq:angmomform} gives
\begin{equation} \label{eq:LangMom}
L=I\dot{\theta}+\tilde{F_1}\dot{\alpha_1}+\tilde{F_2}\dot{\alpha_2}+\tilde{F_3}\dot{\alpha_3},
\end{equation}
where $I$ is the moment of inertia as found in equation \eqref{eq:moialpha} and 
\begin{subequations}
\begin{align}
\tilde{F_1}  & =   2+\cos{\alpha_1}\cos{\alpha_2}+\frac{4}{5}\cos{\alpha_1}\cos{\alpha_3}+\frac{6}{5}\cos{\alpha_2}\cos{\alpha_3}\nonumber\\
& \quad  +\frac{3}{5}\sin{\alpha_1}\sin{\alpha_2}+\frac{6}{5}\sin{\alpha_2}\sin{\alpha_3},
\\
\tilde{F_2}  & =   2+\cos{\alpha_1}\cos{\alpha_2}+\frac{4}{5}\cos{\alpha_1}\cos{\alpha_3}+\frac{6}{5}\cos{\alpha_2}\cos{\alpha_3}\nonumber\\
& \quad  +\frac{3}{5}\sin{\alpha_1}\sin{\alpha_2}+\frac{4}{5}\sin{\alpha_1}\sin{\alpha_3},
\\
\tilde{F_3}  & =   \frac{12}{5}+\frac{8}{5}\cos{\alpha_1}\cos{\alpha_2}+\frac{4}{5}\cos{\alpha_1}\cos{\alpha_3}+\frac{6}{5}\cos{\alpha_2}\cos{\alpha_3}\nonumber\\
& \quad  +\frac{8}{5}\sin{\alpha_1}\sin{\alpha_2}.
\end{align}
\end{subequations}

The same expression for the angular momentum $L$ in terms of orientation and  shape coordinates $(\theta, \alpha_1, \alpha_2, \alpha_3)$ 
can be derived by starting with the Lagrangian given by the kinetic energy $ {\cal L} = \frac12 \sum m_i |\dot z_i|^2$. 
After introducing orientation and shape coordinates the Lagrangian becomes independent of $\theta$, 
and the conjugate momentum $\partial {\cal L} / \partial \dot\theta$ is the conserved angular momentum.
In this derivation we would either treat $\alpha_3$ as a known function given by (14) or use Lagrange multipliers.
In the free motion of the system the shape and orientation would be determined by the corresponding 
Euler-Lagrange equations. 
In the following we instead consider the shape as given by explicitly time-dependent functions $\alpha_i(t)$, 
and the only equation of motion we use is \eqref{eq:LangMom} to find $\theta(t)$.

\section{Geometric Phase}

In the coordinates $(\theta, \alpha_1, \alpha_2, \alpha_3)$ we can interpret the formula for the angular momentum 
as being decomposed into a contribution that originates from changing the shape, 
and a single term $I \dot \theta$ that originates from the rotation of the shape. 
For a rigid body this would be the only term present.
The well known but nevertheless surprising result is that even when $L = 0$ 
the orientation $\theta$ may change. This is most clearly seen when we solve for $\dot \theta$, which gives
\begin{equation}
\dot{\theta}=\frac{L}{I} - \left(\frac{\tilde{F_1}}{I}\dot{\alpha_1}+\frac{\tilde{F_2}}{I}\dot{\alpha_2}+\frac{\tilde{F_3}}{I}\dot{\alpha_3}\right)=
\frac{L}{I} + F_1\dot{\alpha_1}+F_2\dot{\alpha_2}+F_3\dot{\alpha_3}.
\label{eq:thetadot}
\end{equation}
This gives the decomposition of the change of $\theta$ into a dynamic phase proportional to the constant angular momentum $L$, 
and a geometric phase proportional to the shape change, which is proportional to the time derivatives of the angles $\alpha_j$.
Henceforth we will set the angular momentum to zero, $L=0$, so that we can study the question of how to 
maximise orientation change in the absence of angular momentum.

The geometric definition of the rotation angle $\theta$ is not unique, and making a particular choice is called a gauge.
Given a particular shape change $\alpha_i(t)$ the resulting overall rotation $\theta(t)$ at a particular time $t$ depends on the gauge.
To get a gauge invariant quantity we consider $\Delta \theta = \theta(t_1)  - \theta(t_0)$ for closed loops $\gamma$ in shape space, 
that is for $\alpha_i(t_1) = \alpha_i(t_0)$. For a closed loop in shape space it makes sense to subtract 
$\theta$'s to get $\Delta\theta$, since they are computed for the same shape, see for example~\cite{LittlejohnReinsch97} for more information about gauge invariance.

Integrating equation \eqref{eq:thetadot} with $L=0$ over a loop $\gamma $ in shape space gives the overall change in $\theta$,
\begin{equation}\label{eq:deltatheta}
	\Delta \theta =\oint \limits_\gamma \! F_1 d\alpha_1 + F_2 d\alpha_2+ F_3 d\alpha_3.
 \end{equation}
It is not obvious how to choose the loop $\gamma$ such that this integral becomes large. 
In order to find the optimal loop we use Stokes' theorem to change the line integral \eqref{eq:deltatheta} into a surface integral
\begin{eqnarray}
\Delta \theta &=& \iint \limits_{S(\gamma)} \! \nabla \times \mathbf{F} \cdot \,d \mathbf{S} \label{eqn:stokes}
\end{eqnarray}
where $\mathbf{F} = (F_1, F_2, F_3)^t$, $\nabla \times$ is the curl with respect to the angles $\alpha_i$, 
and $S(\gamma)$ is the surface enclosed by $\gamma$ on shape space $C = 0$. 
The formula for $\Delta \theta$ is an integral over a two-form (after eliminating, say, $\alpha_3$). 
We want to convert this integral over a two-form into an integral over a function on shape space.
The function thus defined will be called the magnetic field $B$. 
In order to do this we need a metric $g$ on our shape space surface, so that the magnetic field becomes invariantly 
defined after dividing the two-form by the area form $\sqrt{ \det g}$.

In general a ``magnetic field'' arises from  reduction by a continuous symmetry.
Like the equation for the angular momentum $L$, the magnetic field could be derived in the Lagrangian or Hamiltonian formalism, 
see for example \cite{MMR90,LittlejohnReinsch97}. 
There, one would use e.g.\ $(\alpha_1, \alpha_2)$ as local coordinates and obtain a magnetic field 
$( \nabla \times {\bf F}) \cdot \nabla C/  C_3$
where  $C_3 = \partial C / \partial {\alpha_3}$ is the derivative of the constraint \eqref{eq:conalpha}.
The natural metric $g$ in our case is obtained from the kinetic energy of  the Lagrangian of the Pentagon {\em after} symmetry reduction, 
so that $\tilde B =  ( \nabla \times {\bf F}) \cdot \nabla C/ ( C_3 \sqrt{ \det g})$ is the magnetic field.
However, the expression for $\det g$ is fairly complicated, and since we are not interested in the free motion of the pentagon
(which is given by the geodesic flow of this metric) we prefer to use a different approach.
Suffice it to say that the magnetic field  $\tilde B$ is invariant
under our discrete symmetry group $D_{10}$ and has zero average over the whole shape space.
The geometric interpretation of this average as a Chern class confirms that the bundle of orientation $\theta$ 
over shape space is a trivial bundle.

Instead we attempt to use the metric induced by the embedding of shape space in the three-dimensional torus $\mathbb{T}^3$ with 
coordinates $\alpha_i$ given by the equation $C=0$, see  \eqref{eq:conalpha}. 
The area form of this surface using $\alpha_1$ and $\alpha_2$ 
as local coordinates is $d S = |\nabla C|/C_3 \, d\alpha_1 d\alpha_2$, 
and with $d {\bf S} = {\bf n} d S$ where ${\bf n} = \nabla C / |\nabla C|$ the integrand in \eqref{eqn:stokes} becomes 
$( \nabla \times {\bf F}) \cdot \nabla C / C_3 \, d\alpha_1 d\alpha_2 = ( \nabla \times {\bf F}) \cdot \nabla C/ | \nabla C| dS$,
so that the magnetic field would be the factor multiplying $dS$.
This function differs from $\tilde B$ by a scalar non-zero factor.
The trouble with the function $( \nabla \times {\bf F}) \cdot \nabla C/ | \nabla C| $ is that it is not invariant under the discrete symmetry group.
This is not too surprising, since we arbitrarily invented a metric on shape space by embedding it into $\mathbb{T}^3$ using $C=0$.
In the following all we want to do with the magnetic field is to study its zero level curve $B=0$ in the fundamental domain $\phi$.
Thus multiplying $\tilde B$ by a non-zero factor does not make a difference. 
In order to make the result invariant under the discrete symmetry group we multiply the function obtained from 
the metric of the surface $C=0$ by $|\nabla C|$, and thus define our magnetic field to be given
by the simple expression
\begin{equation}
   B = ( \nabla \times \mathbf{F}) \cdot \nabla C \,.
\end{equation}
Explicitly we find 
\[
\begin{split}
\frac{5}{8} I^2 B = 
\cos (3 \alpha_3) \big[ \cos \alpha_1+\cos \alpha_2\big] +\cos (2 \alpha_3) \big[3 + 2 \cos( \alpha_1 - \alpha_2) + \cos( \alpha_1 + \alpha_2)\big]+ \\
+\cos \alpha_3 \big[\cos \alpha_1+\cos \alpha_2-\cos (3 \alpha_2)-\cos (2 \alpha_1-\alpha_2)-2 \cos (\alpha_1+2 \alpha_2) \big]+\\
   +\sin \alpha_1 \sin \alpha_2-\cos \alpha_1 (2 \cos \alpha_2+\cos (3 \alpha_2))-\cos (2 \alpha_1)-2 \cos (2 \alpha_2) \,.
\end{split}
\]
The  ``magnetic field'' $B$ is invariant under the the action of the discrete symmetry group $D_{10}$.
We don't have a good explanation why this expression possesses the correct discrete symmetry, 
but this together with its simplicity if the main reason to use it instead of $\tilde B$.
Using the proper magnetic field $\tilde B$ instead of $B$ does not change our results, 
only a different parametrisation of the same loop $B = 0$ would be used.


\section{Optimal Shape Change}

We are seeking the optimal contractible loop $\gamma$ on shape space in the sense that it maximises overall orientation change 
of the pentagon given by \eqref{eq:deltatheta} after one revolution. 
Imagine we start with a small loop $\gamma$ in a region where $B > 0$, say near a positive maximum of $B$.
Enlarging the loop will increase $|\Delta \theta|$ as long as $B > 0$. 
This process can be repeated and we can keep growing the loop, yielding a larger and larger $|\Delta \theta|$. 
However, enlarging the loop $\gamma$ across the $B = 0$ contour would give an opposing contribution towards the integral, thus lowering $|\Delta \theta|$.
Hence the largest $|\Delta \theta|$ is achieved when the loop $\gamma$ coincides with the contour $B = 0$, assuming there is a contractible 
zero-contour of $B$ enclosing the initially chosen small loop.
In general, the connected components of the $B=0$ contour on a genus 4 surface may be non-contractible. 
Yet in the present case all zero-contours of $B$ are contractible.

The sign of $\Delta \theta$ is dependent on the orientation of $\gamma$. A positively oriented loop in the $B>0$ region will yield a negative $\Delta \theta$. Though this sounds counter intuitive, the reason is that the matrix $A$ from the affine transformation given in \eqref{eq:psi2alpha} has a negative determinant. That is, a positively oriented loop in $\alpha$\textit{-space} corresponds to a negatively oriented loop in $\psi$\textit{-space}. The construction works similarly starting in an initial region with $B<0$, where a positively oriented loop will yield a positive $\Delta \theta$. 

The connected components of the zero-contours of $B$ give analytic curves that yield locally optimal loops 
that make $\Delta\theta$ extremal, in the sense that any small variation of the loop  decreases the value of $| \Delta \theta|$. 
Globally it is possible to connect two disjoint regions with $B\ge0$ with a curve through the area with $B < 0$. 
Traversing this curve back and forth provides no net contribution as it does not enclose any area.
Alternatively one could also traverse the same loop twice and hence double the amount of rotation.

\begin{figure}[tb]
\centering
\includegraphics[width=9cm]{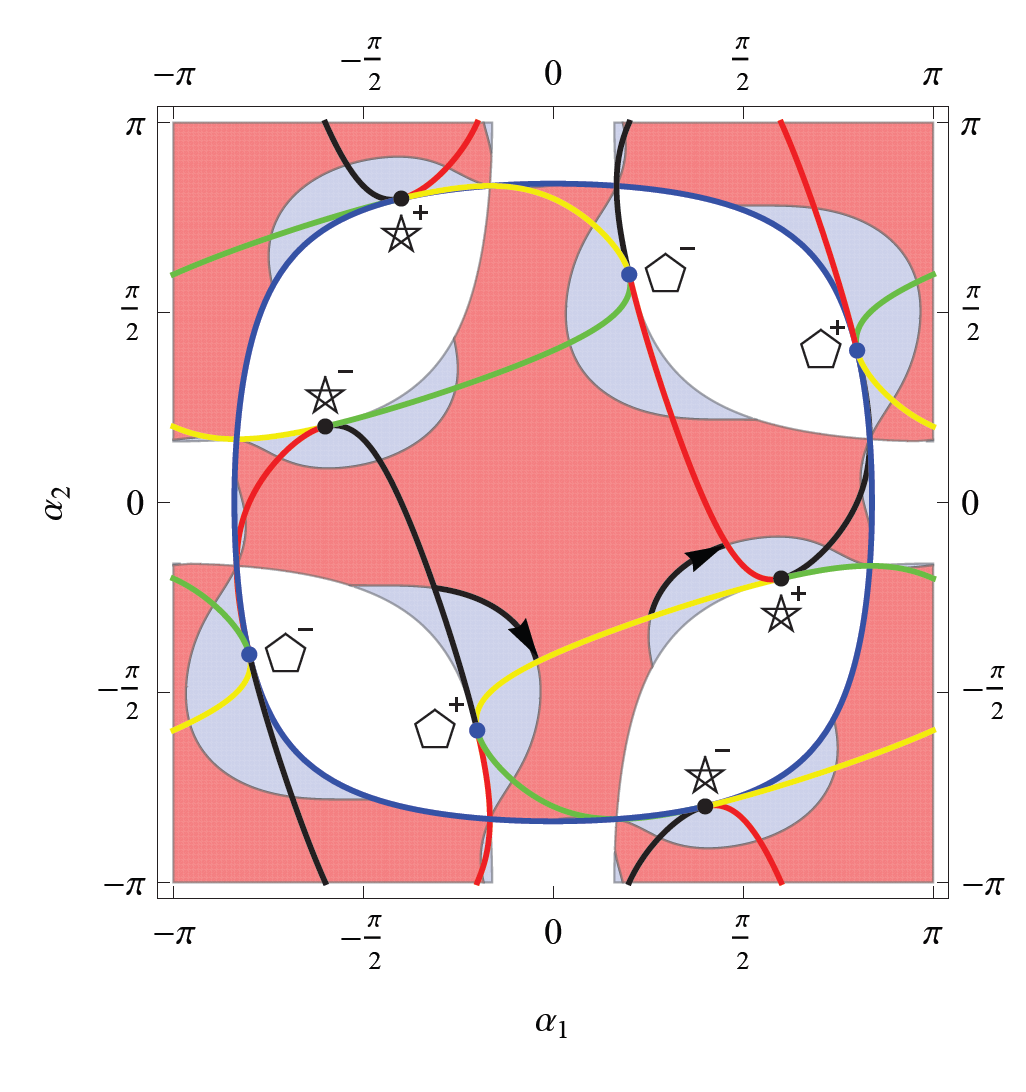}
\caption{The light blue region indicates $B>0$ and pink $B<0$. 
There are 4 curves with $B=0$ in the fundamental region, enclosing the positively/negatively oriented regular pentagon/pentagram, respectively.
\newline An interactive Mathematica file {\em link to manipulate.cdf} allows you to select a point in the fundamental region $\phi$ 
with the mouse cursor and displays the corresponding pentagon.}\label{fig:bzero}
\end{figure}

It is enough to consider the $B = 0$ contour in the fundamental region $\phi$ because of the discrete symmetry, as shown in Figure \ref{fig:bzero}.
The symmetry $V$ maps curves enclosing the positively oriented regular pentagon/pentagram into curves 
enclosing the negatively oriented regular pentagon/pentagram, respectively. As $V$ preserves $\Delta \theta$, only the two loops shown in Figure~\ref{fig:bzero} need to be considered. The orientation is as indicated and specifically chosen such that the overall rotation $\Delta \theta$ is positive.

When defining the optimal loop $\gamma$ through $B=0$ we may multiply $B$ by any smooth function on the shape space
that is invariant under the symmetry group $D_{10}$ without changing the result.
The overall denominator of $B$ is the moment of inertia squared, which can be removed when computing $\gamma$.
The integrand in \eqref{eqn:stokes} involves rational functions of trigonometric functions of $\alpha_i$ 
that is even in $\alpha_3$. Thus we can replace $\alpha_3$ by $\alpha_3( \alpha_1, \alpha_2)$ 
and obtain a rational function of trigonometric functions of $\alpha_1, \alpha_2$. 
We have not been able to find a simple closed form for this integral, so we resort to numerical methods at this point.

\begin{figure}[tb]
\centering
\includegraphics[width=14cm]{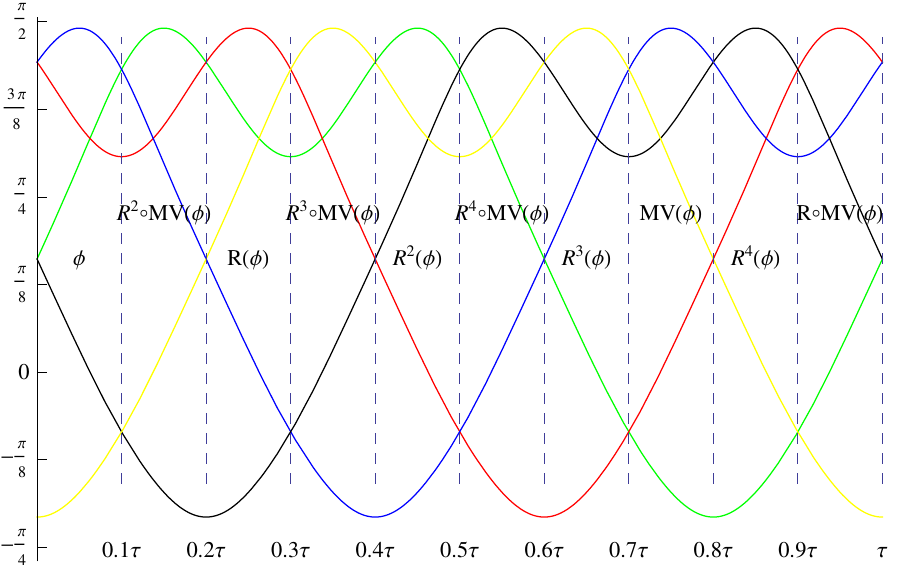}
\caption{The diagram shows how the relative angles $\psi_i$ change along the $B=0$ contour around the pentagram, the colouring of $\psi_i$ is from our standard colour code, that is, $i=1,\dots, 5$ is green, red, blue, black, yellow. 
The curve segment in $\phi$ as seen in Figure \ref{fig:bzero} was used to construct the entire loop on shape space using the indicated symmetry operations.}\label{fig:symcon} 
\end{figure}

Once the curve is known in the fundamental region the symmetries are used to find the curve in all of shape space. 
The result of this is shown in Figure~\ref{fig:symcon}.
In order to parameterise the zero-contour of $B$, we numerically solve Hamilton's equations with a
Hamiltonian given by $H = B I^2$ with initial conditions on the boundary of the fundamental region $\phi$.
Note that this Hamiltonian is merely used for the purposes of obtaining a parametrisation of the $B=0$ contour; 
it is not the Hamiltonian of the free motion of the pentagon.
The speed with which the contour is traversed is irrelevant for the final geometric phase $\Delta\theta$,
so whether we consider $H = B$ or $H = B I^2$ makes no difference. We can even multiply by $B$ by a
function that is not invariant under the $D_{10}$ action, 
but then the overall solution pieced together from the action of the symmetry group 
on the fundamental piece may not be smooth but only once differentiable.
The ODE in local coordinates $(\alpha_1, \alpha_2)$ is 
\[
    \dot \alpha_1 =  -S \frac{ \partial H }{ \partial \alpha_2}, \quad
    \dot \alpha_2 =   S \frac{ \partial H }{ \partial \alpha_1} \,.
\]
Here $H = H(\alpha_1, \alpha_2, \alpha_3(\alpha_1, \alpha_2))$ and $S = C_3(\alpha_1, \alpha_2, \alpha_3(\alpha_1, \alpha_2))$ 
is the symplectic multiplier from the area-element. 
Maybe the easiest way to derive this non-standard symplectic structure is to start with the global ODE on the three-torus
$\dot \alpha = \nabla C \times \nabla B$ (which also gives an alternative way to numerically compute the whole loop $B=0$)
and consider $\nabla C \times$ as a Poisson structure. 
Reduction to the symplectic leaf $\{ C = 0 \}$ using local coordinates $(\alpha_1, \alpha_2)$ then gives the symplectic structure 
$S d\alpha_1 \wedge d\alpha_2$. 

The solution curves $\alpha_i(t)$, and hence $\psi_i(t)$ as well, along the $B=0$ contour are periodic with period $\tau$.
The solutions for the different angles are related by a phase shift, thus it is enough to study a single curve for the whole period $\tau$, 
say $\psi_5(t)$, which is an even function.
Specifically the positively oriented loop in $\alpha$-space around the positive pentagram has the relation $\psi_i(t) = \psi_{i+1}( t + \frac{\tau}{5})$ while the negatively oriented loop has relation $\psi_i(t) = \psi_{i+1}( t - \frac{\tau}{5})$. The phase shift relation for the positively and negatively oriented loops in $\alpha$-space around the positive regular convex pentagon are $\psi_i(t) = \psi_{i+1}( t + \frac{2\tau}{5})$ and $\psi_i(t) = \psi_{i+1}( t - \frac{2\tau}{5})$, respectively.

A natural way to encode the final answer is to decompose $\psi_5(t)$ into a Fourier cosine series,
\begin{equation}
\psi_5(t)=\displaystyle\sum_{n=0}^{\infty}{\tilde a_n \cos{(n \omega t)}} \mbox{ \indent where } \omega = \frac{2\pi}{\tau} = 7.3634...,
\end{equation}
and the coefficients $\tilde a_n$ are given in Table \ref{table:fourieran}

\begin{table}[ht]
\centering 
\begin{tabular}{c c c c} 
\hline\hline 
$n$ & $\tilde a_n$ & $n$ & $\tilde a_n$\\ [0.5ex] 
\hline 
0 & 0.6283 & & \\
1 & -0.9646 & 8 & -0.004078\\
2 & -0.3974 & 9 & 0.0003710\\
3 & 0.1595 & 10 & 0\\
4 & -0.0779 & 11 & -0.0002426\\
5 & 0 & 12 & 0.0001005\\
6 & 0.008748 & 13 & 0.00003240\\
7 & -0.002677 & 14 & 0.00001171\\[1ex] 
\hline 
\end{tabular}
\caption{Numerical values of the Fourier coefficients $\tilde a_n$ of $\psi_5(t)$. Along this loop where $B=0$ truncating the Fourier series at $n=14$ produces an absolute error $|B(\psi_i(t))|$ on the order of $10^{-5}$.} 
\label{table:fourieran} 
\end{table}

The fact that all coefficients $\tilde a_{5n} = 0$ is equivalent to $\sum_{j=0}^{4} \psi_i(t + \frac{j \tau}{5}) = const = 5 \tilde a_0$.
Using the property $\psi_i(t) = \psi_{i+1}(t + \frac{\tau}{5})$ translates this into $\sum \psi_i(t) = 5 \tilde a_0$.
By construction we have that $\sum \psi_i = \pi \bmod 2 \pi$, hence $\tilde a_0 = \pi/5$.
If we consider the sub-Fourier series defined by $f_l(t) = \sum a_{l + 5 n } \cos( (l + 5 n ) \omega t)$ so that
$f_0(t) = \frac{\pi}{5}$, we find that all angles can be explicitly written as linear combinations of four functions $f_l$, $l = 1, 2, 3, 4$
(and the constant function $f_0$) as
\[
    \psi_{5 + j}(t) = \Re \sum_{l=0}^{4} f_l(t) e^{2\pi i j l/5} \,.
\]
The function $f_l(t)$ can be written as $e^{i \omega l t} \hat f_l(t)$ where $\hat f_l(t) = \sum c_{l+kn}e^{i\omega k n t}$ has period $\frac{\tau}{5}$.
In real variables this becomes 
$\hat f_l(t) = \cos( l \omega t) ( a_{l, 0} + a_{l, 1} \cos( 5 \omega t) + \dots)$  where $a_{l,n} = a_{l + 5 n }$.
Any periodic shape change that is obtained by unfolding a curve from the fundamental region has these properties. 


As $\dot\theta$ is a composition of periodic functions, 
it is also a periodic function with the same period $\tau$, with Fourier series 
\begin{equation}\label{eq:thetadotfourier}
\dot{\theta}(t) = a_0+\displaystyle\sum_{n=1}^{\infty} {\left\{a_n \cos{(n\omega t)}+b_n \sin{(n \omega t)}\right\}},
\end{equation}
where the numerical values of $a_n$ and $b_n$ are given in Table \ref{table:fourierthetatable}.
Notice that even for the optimal loop the sign of $\dot \theta$ is not constant along the loop. 
\begin{table}[ht]
\centering 
\begin{tabular}{c c c} 
\hline\hline 
$n$ & $a_n$ & $b_n$\\ [0.5ex] 
\hline 
0 & 0.9239 & -\\ 
1 & -1.8671 & -5.7463\\
2 & -2.4891 & 1.8085\\
3 & -1.4990 & -1.0891\\
4 & 0.6032 & -1.8566\\
5 & -0.3977 & 0\\
6 & 0.1016 & 0.3127\\
7 & -0.0587 & 0.04265\\
8 & 0.1022 & 0.07423\\
9 & -0.006463 & 0.01989\\
10 & 0.005644 & 0\\
[1ex] 
\hline 
\end{tabular}
\caption{The Numerical Values of the $a_n$ and $b_n$ terms of $\dot{\theta}$} 
\label{table:fourierthetatable} 
\end{table}

\noindent
Integrating equation \eqref{eq:thetadotfourier} term by term gives
\begin{equation}\label{eq:thetafourier}
\theta (t)=C+a_0 t+\displaystyle\sum_{n=1}^{\infty} {\left\{\frac{a_n}{n\omega} \sin{(n \omega t)}-\frac{b_n}{n \omega} \cos{(n \omega t)}\right\}}.
\end{equation}
Hence $\Delta\theta = \theta(\tau) - \theta(0) = a_0 \tau$.
The numerical value for $\Delta\theta$ obtained along the locally optimal loop around the regular pentagram is
\begin{equation}
\Delta\theta=a_0 \tau \approx 0.78837 \approx 45.17^{\circ} \,.
\end{equation}
Performing the same analysis for the locally optimal loop around the regular convex pentagon yields
\begin{equation}
\Delta\theta \approx 0.49147 \approx 28.16^{\circ}.
\end{equation}

\begin{figure}[tp]
\centering
\subfigure[t=0]{
\includegraphics[width=2cm]{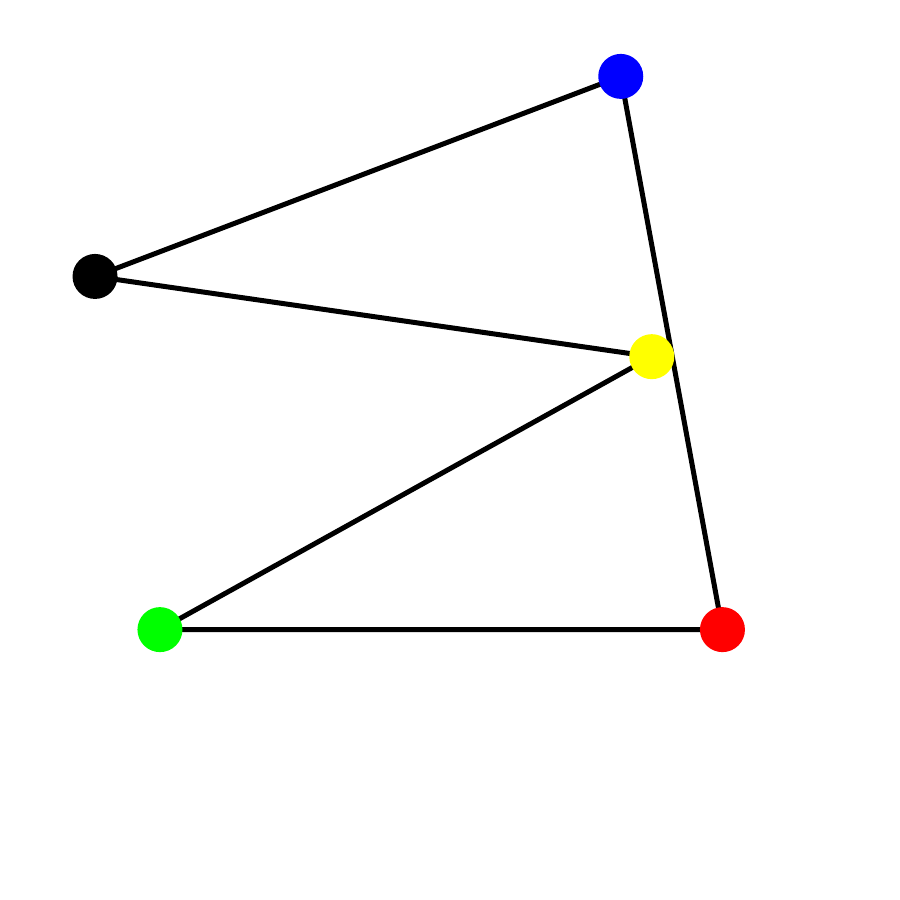}
}
\subfigure[$t=\frac{\tau}{100}$]{
\includegraphics[width=2cm]{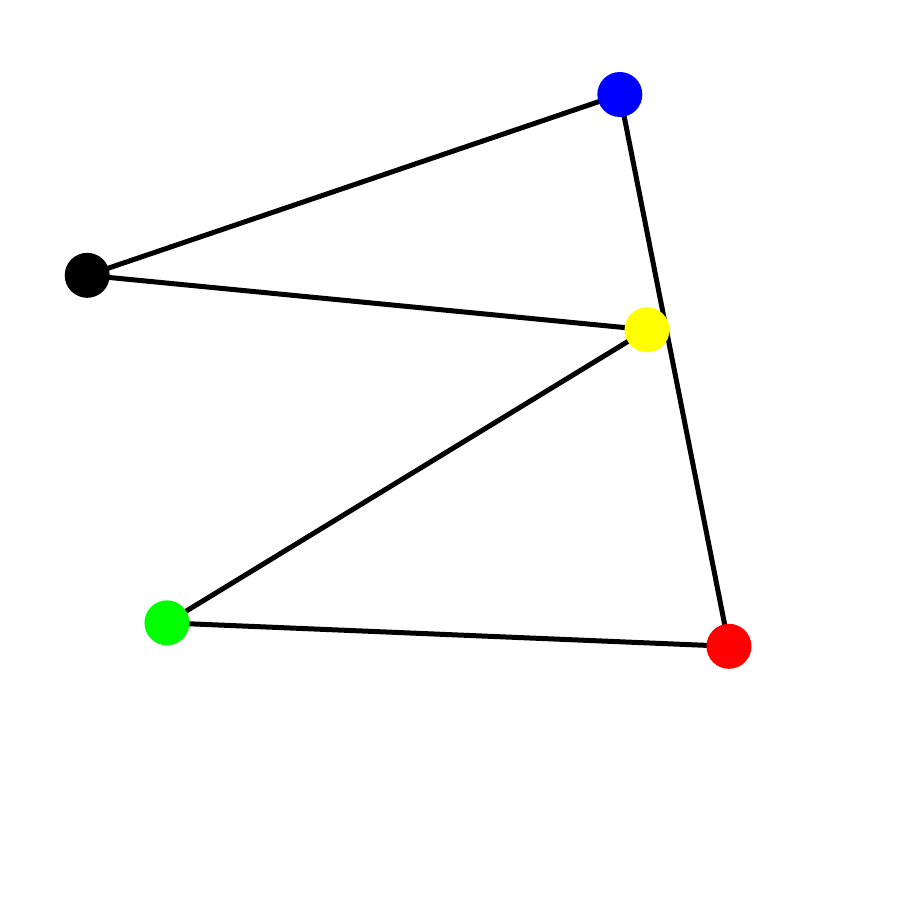}
}
\subfigure[$t=\frac{2\tau}{100}$]{
\includegraphics[width=2cm]{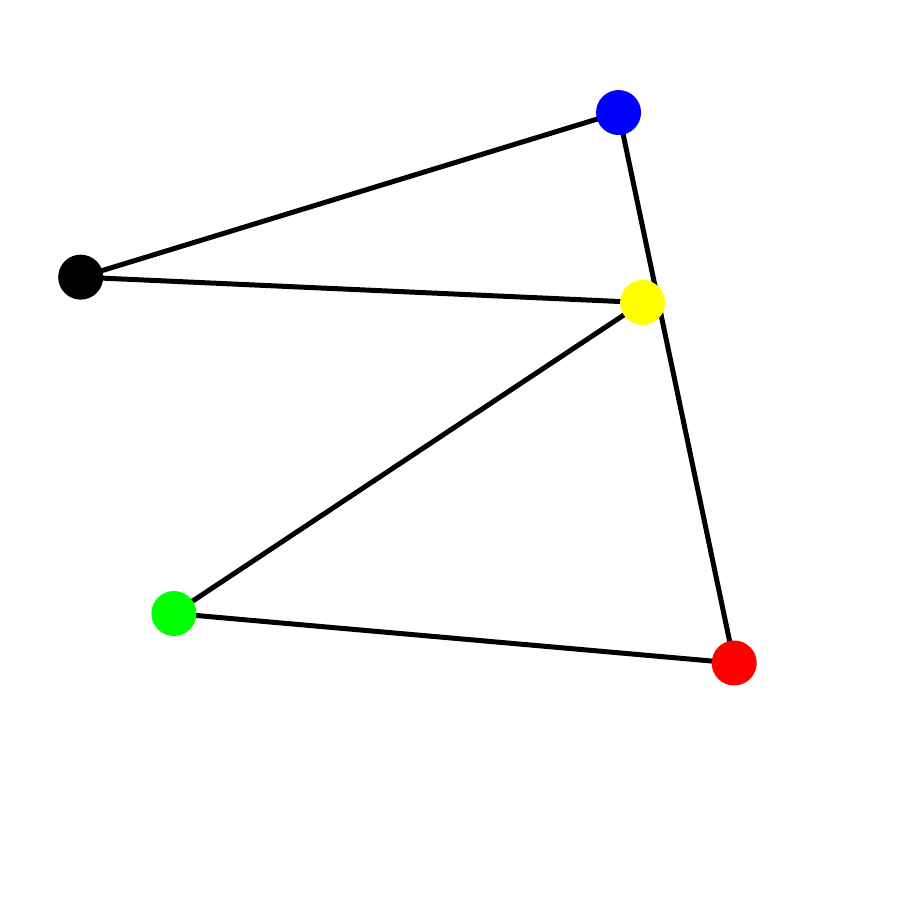}
}
\subfigure[$t=\frac{3\tau}{100}$]{
\includegraphics[width=2cm]{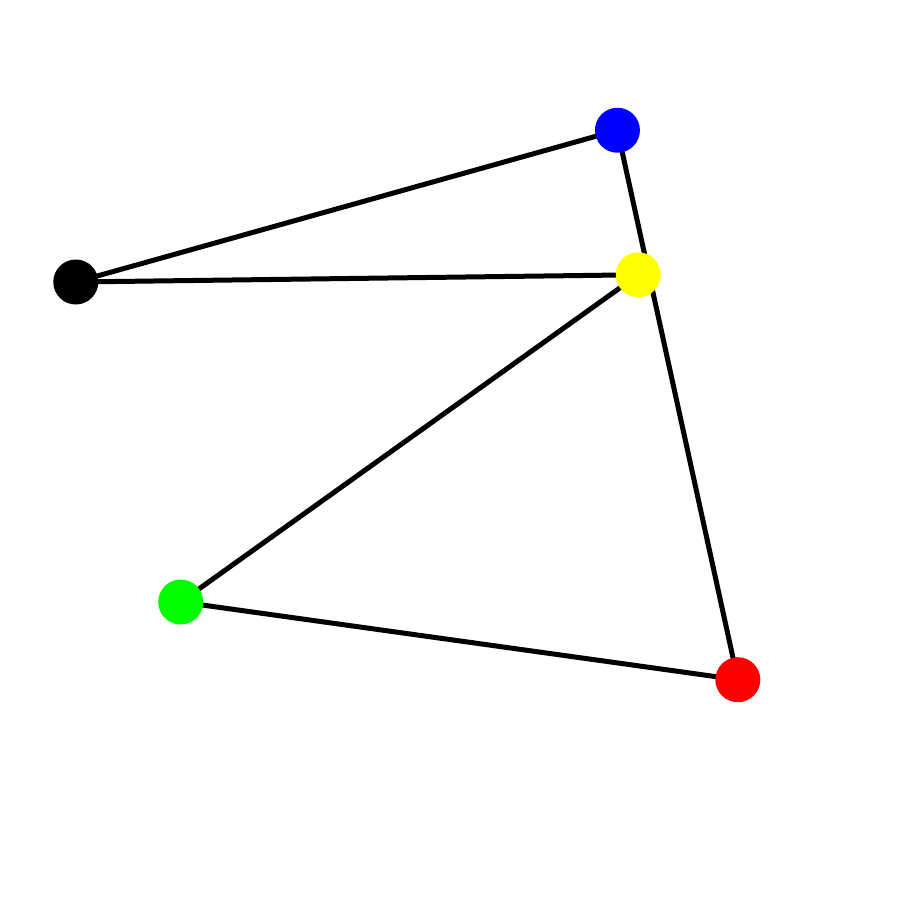}
}
\subfigure[$t=\frac{4\tau}{100}$]{
\includegraphics[width=2cm]{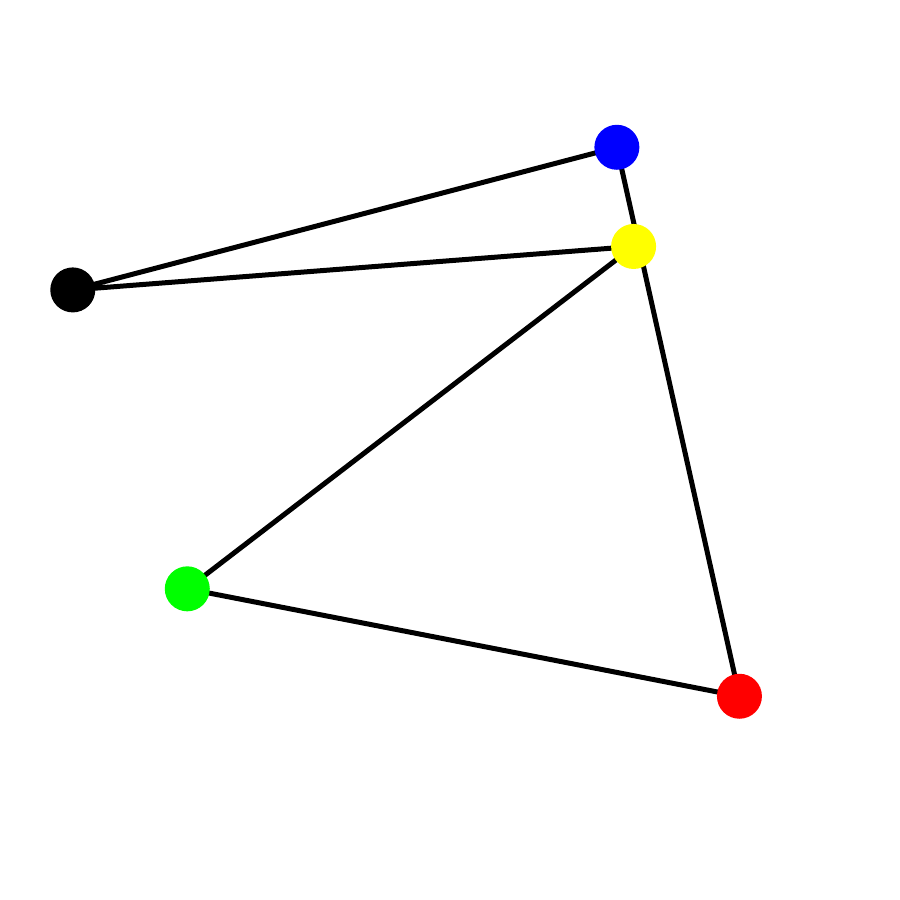}
}
\subfigure[$t=\frac{5\tau}{100}$]{
\includegraphics[width=2cm]{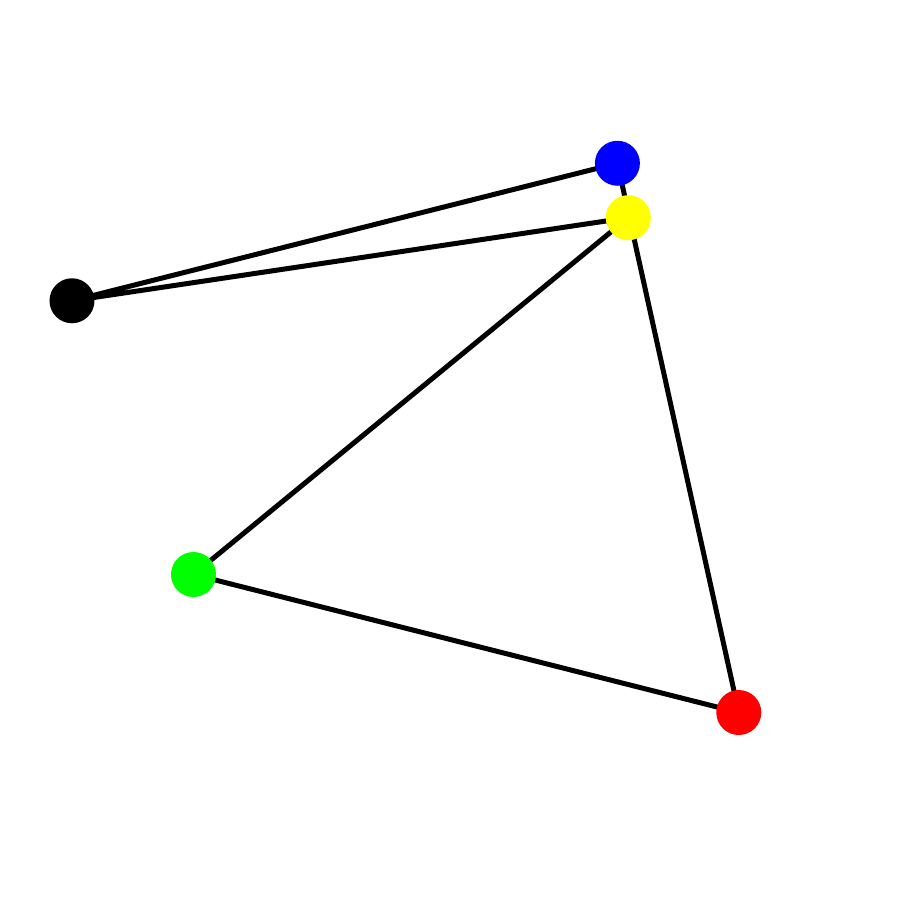}
}\\
\subfigure[$t=\frac{6\tau}{100}$]{
\includegraphics[width=2cm]{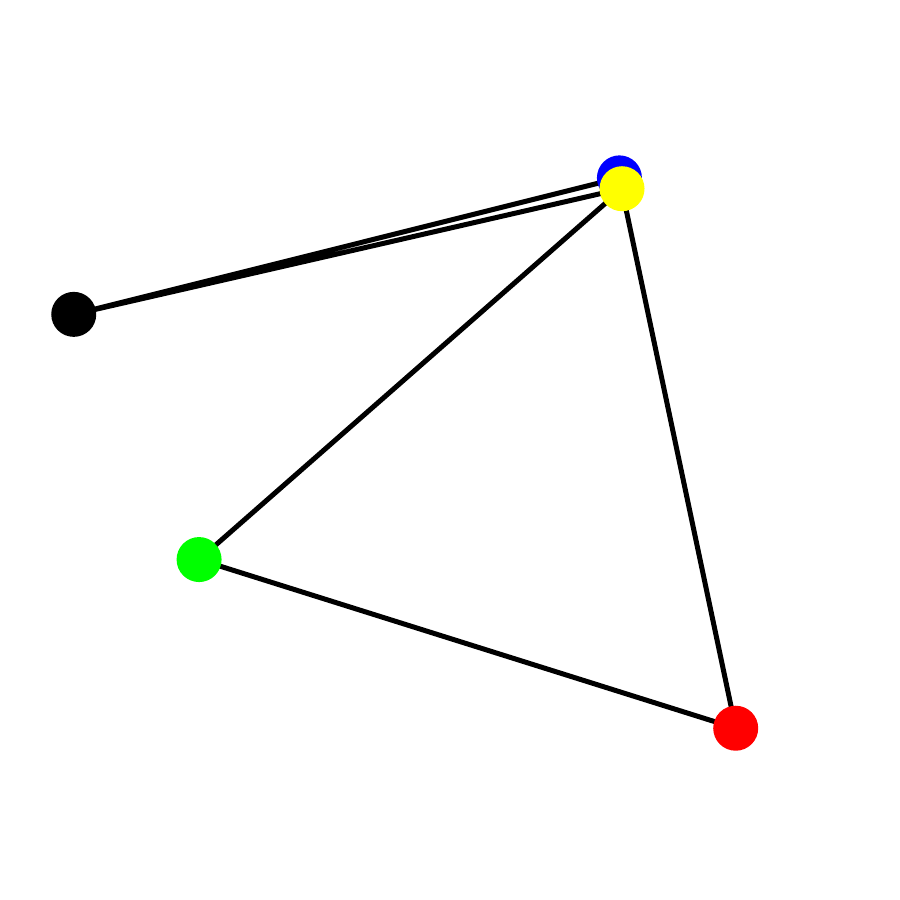}
}
\subfigure[$t=\frac{7\tau}{100}$]{
\includegraphics[width=2cm]{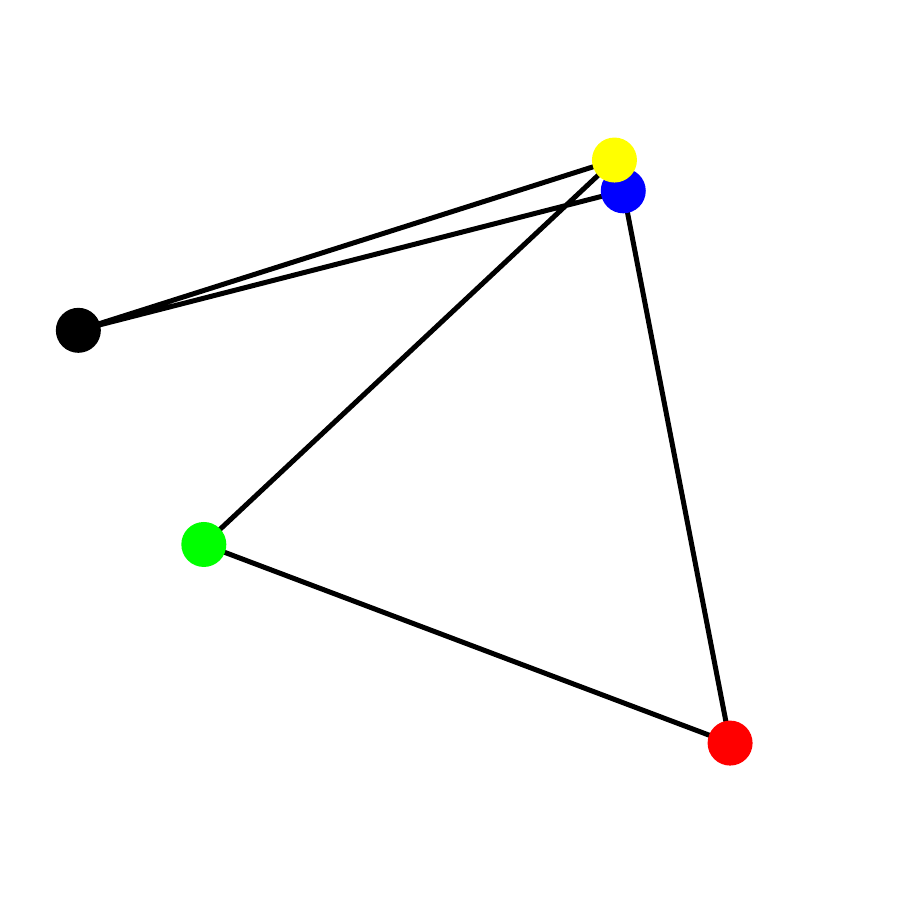}
}
\subfigure[$t=\frac{8\tau}{100}$]{
\includegraphics[width=2cm]{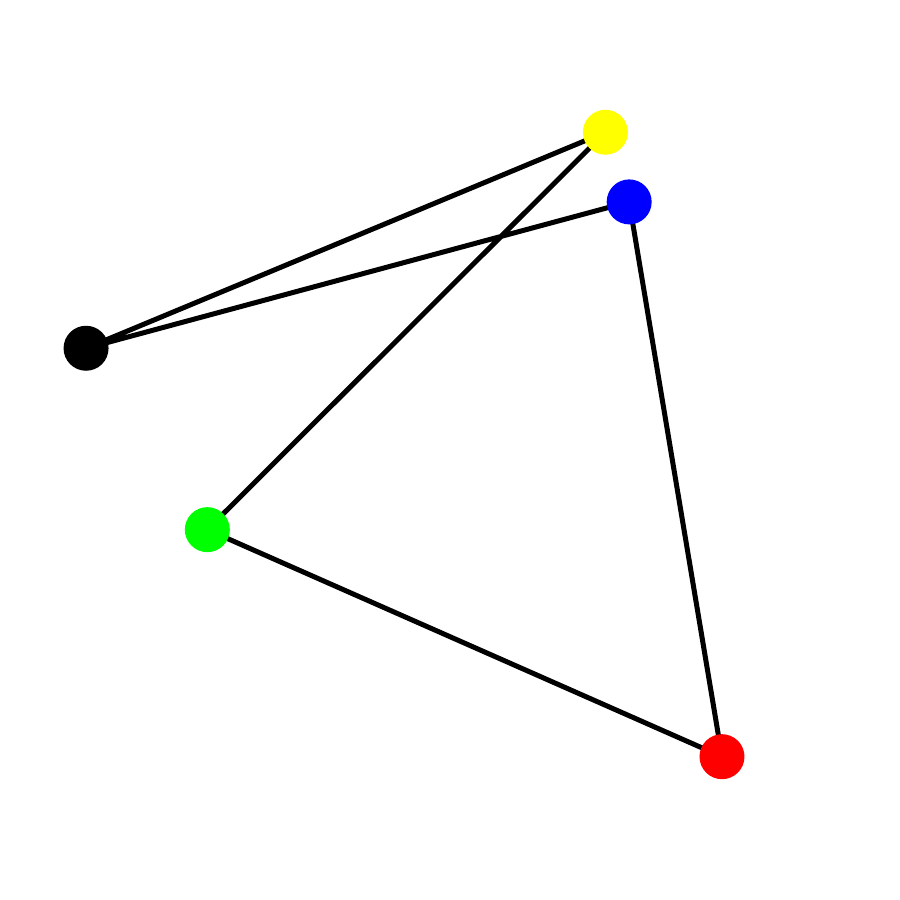}
}
\subfigure[$t=\frac{9\tau}{100}$]{
\includegraphics[width=2cm]{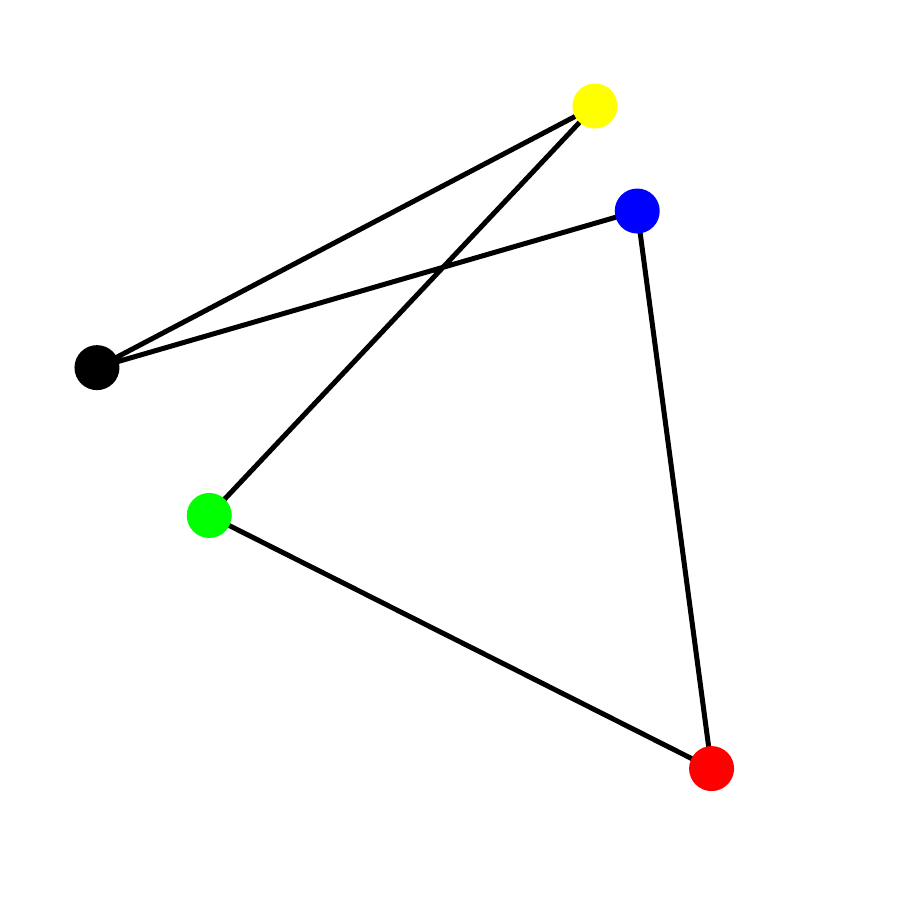}
}
\subfigure[$t=\frac{10\tau}{100}$]{
\includegraphics[width=2cm]{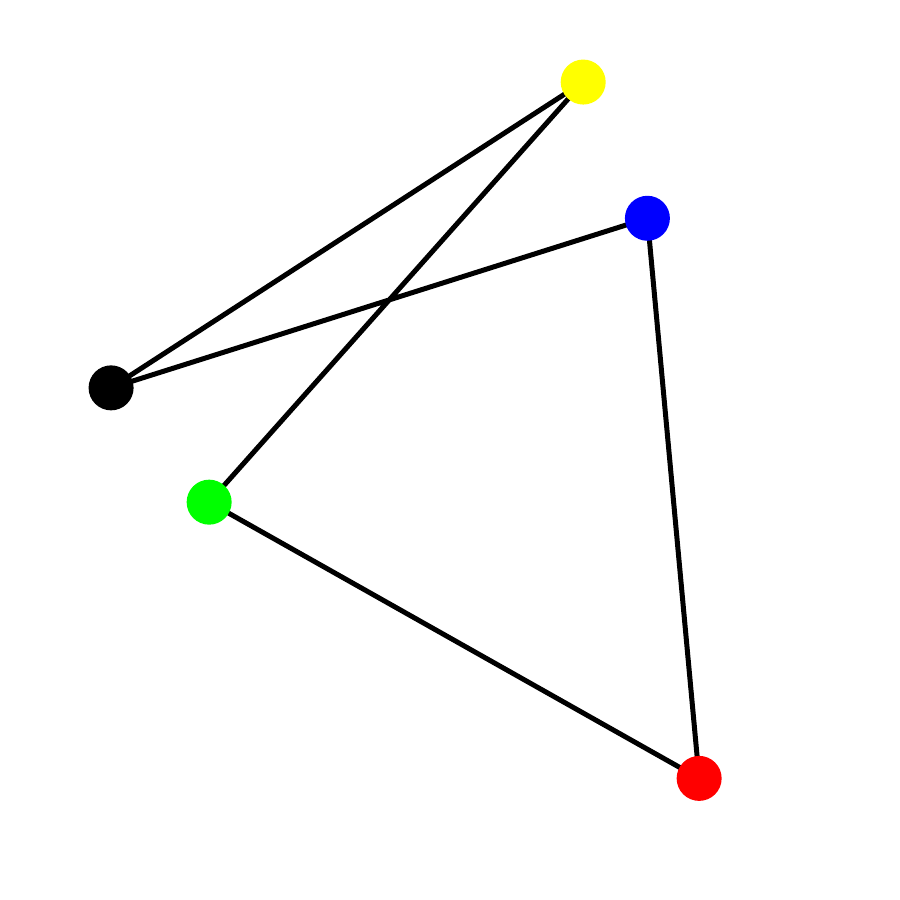}
}
\subfigure[$t=\tau$]{
\includegraphics[width=2cm]{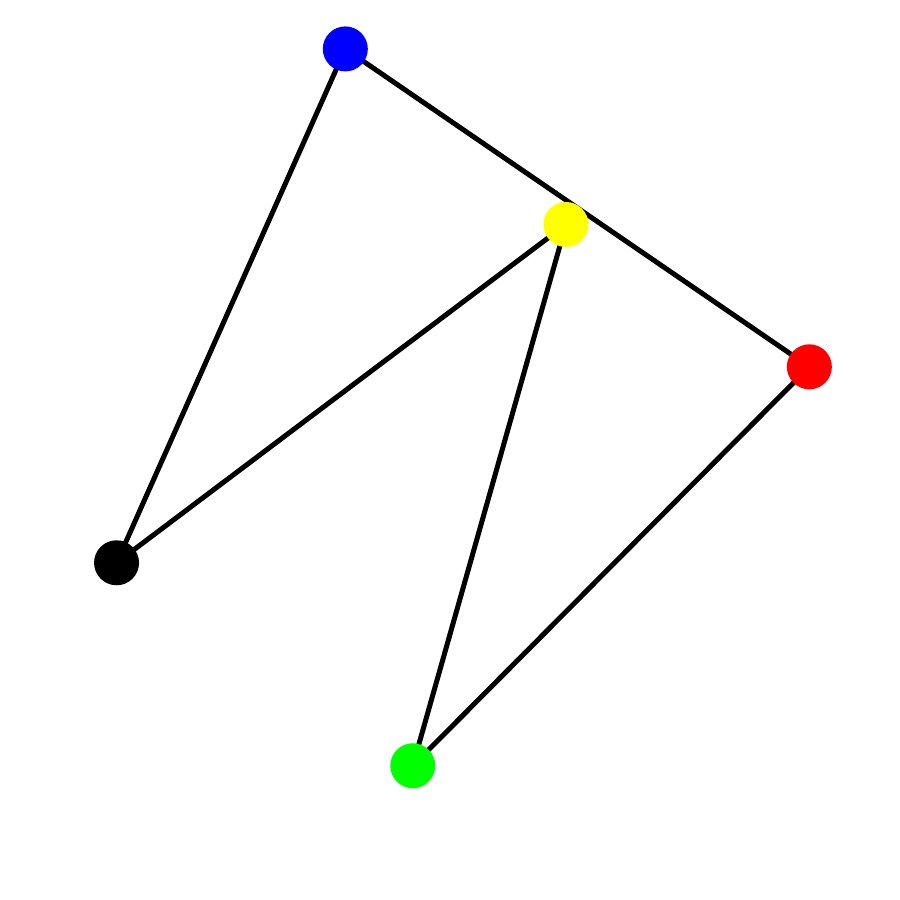}
}
\caption{
The optimal shape change of the equilateral pentagon that maximises overall rotation at angular momentum zero after one period. The diagram shows the shape change within the fundamental region $\phi$ only (except for the last image which serves for comparison between the initial and final orientations of the pentagon). By symmetry reconstruction, the remaining shape changes outside the region continues on with this motion in the reverse order; then the initial motion shown is repeated, each time with different colouring of the vertices until it returns to its original state. 
\newline Movie files {\em links} and/or the Mathematica {\em link to animations.cdf} file shows the complete animations along the different optimal loops of the equilateral pentagon.
}\label{fig:time}
\end{figure}

We conclude this section by illustrating how the shape of the pentagon changes as the loop of optimal shape change is traversed, see Figure \ref{fig:time}.
Although it may appear that the yellow vertex travels along the blue-red edge this is only approximately true.

\section{Convex Pentagons}
The results obtained so far apply to the family of all equilateral pentagons including degenerate and non-simple pentagons. We now restrict to the sub-manifold of equilateral \emph{convex} pentagons.
The intuitive approach to finding the sub-manifold containing all equilateral convex pentagons from the original manifold is to consider the boundary cases for which the pentagon is not strictly convex.
Fixing one of the angles of the equilateral pentagon to $\pi$  gives the curve on shape space which contains all possible equilateral pentagons with the chosen angle $\psi_k=\pi$ for some $k \in \Z_5$. The importance of this curve is that parts of it separate the strictly convex pentagons from the nonconvex pentagons locally.
Choosing $\psi_3=\pi$ we can see from \eqref{eq:psi2alpha} that $\alpha_1=\alpha_2$, so one way to parametrise this curve is to set $\alpha_3=-t$, substitute in \eqref{eq:conalpha} and solving for $\alpha_1=\alpha_2$ gives 
\begin{equation}
p(t)=\left\{\arccos{\chi(t)},\arccos{\chi(t)},-t \right\} \mbox{ where }\chi(t)=-\cos{t}+\sqrt{\cos^2{t}-\frac{3}{4}}.
\label{eq:pcur}
\end{equation}
To orient the curve in the clockwise direction as $t$ increase, we included a minus sign in the third component of the parametrisation in \eqref{eq:pcur}, the first two components are even functions so the minus sign may be omitted.

To find the endpoints of the interval where \eqref{eq:pcur} gives the boundary separating convex and non-convex pentagons, we need another angle $\psi_j$, $j \neq 3$ to equal $\pi$. This occurs when $t=\pm \arccos{\left(\frac{3\sqrt{3}}{4\sqrt{2}}\right)}=\pm p_s$, thus $t\in\left[-p_s,p_s\right]$ is the interval of interest. The pentagonal shapes at these endpoints are the isosceles triangles with side lengths \{1,2,2\}. 

\begin{figure}[h!t]
\centering
\subfigure[]
{\includegraphics[width=7cm]{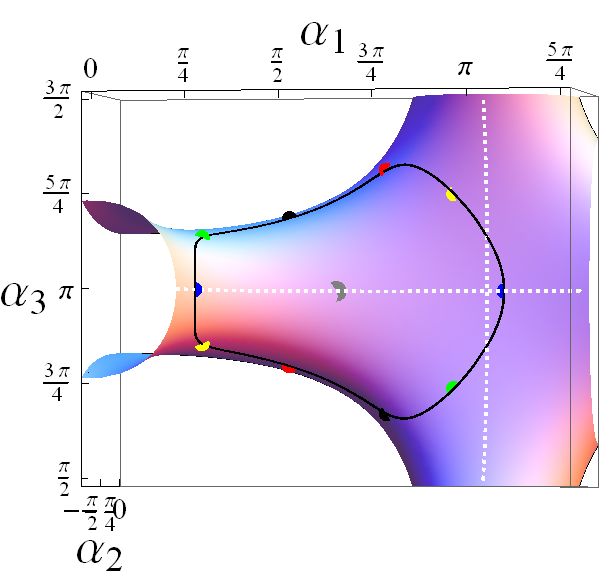}
\label{subfig:starcurve}
}
\subfigure[]
{\includegraphics[width=7cm]{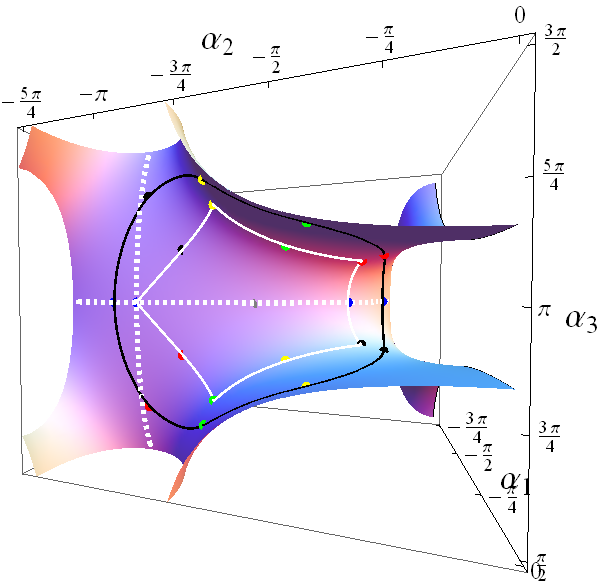}
\label{subfig:pentcurve}
}
\caption{Two views of shape space $C=0$ showing the complete optimal loops.
The gray point in the interior of the loops denotes the pentagram (left) and regular convex pentagon (right). 
Each loop intersects the boundaries of copies of the fundamental region at 10 different points denoted by the coloured points, the colour is chosen to match the corresponding symmetry curve, compare Figure~\ref{fig:map}. The dashed white lines denotes the intersection with $\alpha_i = \pm \pi$.
}\label{fig:views}
\end{figure}

\begin{SCfigure}
\centering
\includegraphics[width=0.41\textwidth]{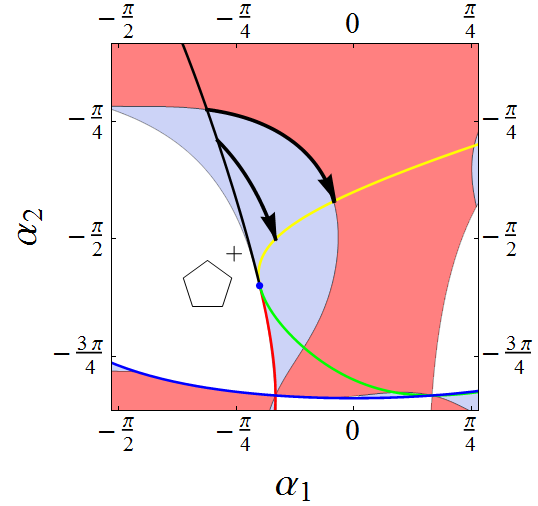}
\caption{This is an enlargement of Figure~\ref{fig:bzero}. In addition the oriented black curve inside the blue region indicates the boundary of convexity given by $R^3(p(t))$. The symmetry operator $R^3$ is applied to $p(t)$ to map the curve to the fundamental region $\phi$. \newline \newline \newline}
\label{fig:zoomedin} 
\end{SCfigure}

Note that when $t\in\left[-p_s,p_s\right]$ the boundary curve traverses through two different fundamental regions, surrounded by the symmetry curves $\{b_1, r_2, k_3, y_4\}$ for $t\in\left[-p_s,0\right]$ and $\{b_1,k_2,r_3,g_4\}$ for $t\in\left[0,p_s\right]$. These two fundamental regions are related by the symmetry operations $MV$. 
To complete the white loop as shown in Figure~\ref{subfig:pentcurve} four additional symmetry operations $R^i$, $i\in\{1,2,3,4\}$  need to be applied to $p(t)$. In order to keep the curve in a single fundamental region we restrict $t$ to $\left[-p_s,0\right]$, 
so that $R^3(p(t))$ for $t \in \left[-p_s,0\right]$ gives the black curve inside $\phi$ with $B>0$ as shown in Figure~\ref{fig:zoomedin}.
The shape at $t=0$ is the trapezium with side lengths $\{1,1,1,2\}$.

The orbit of $p(t)$ under $D_5^+$ gives the piecewise smooth loop around the regular convex pentagon shown in white in Figure~\ref{subfig:pentcurve}.
The particular order of symmetry operations can be read off from Figure~\ref{fig:map}. 
The corner points where the loop is non-smooth are the isosceles triangles that have two angles $\psi_j$ equal to $\pi$, 
while the midpoint of the white edges indicated by the coloured points correspond to the trapeziums with side lengths \{1,1,1,2\}. 
The corresponding loop enclosing convex pentagons of the opposite orientation is obtained by applying either $M$ or $V$.
Also shown in Figure~\ref{fig:views} are the $B=0$ contours in black. 
Intersections of the loops with the boundaries of the 10 fundamental regions that meet at the central positively oriented pentagram (left) 
or positively oriented regular convex pentagon (right) are marked by coloured dots, also compare Figure~\ref{fig:map}.

The white loop which demarcates the boundary of the region of convex pentagons lies completely inside the region 
with $B > 0$ (given by the black loop). Thus the maximal rotation of the equilateral {\em convex} pentagons 
is obtained along this loop. The argument is similar to that used to show optimality for the $B=0$ loop.
Deforming the loop into the non-convex region is not allowed, while deforming the loop in the other direction
would decrease the enclosed area, and hence the overall amount of rotation as $B > 0$.

\begin{figure}[tp]
\centering
\subfigure[]
{\includegraphics[width=7cm]{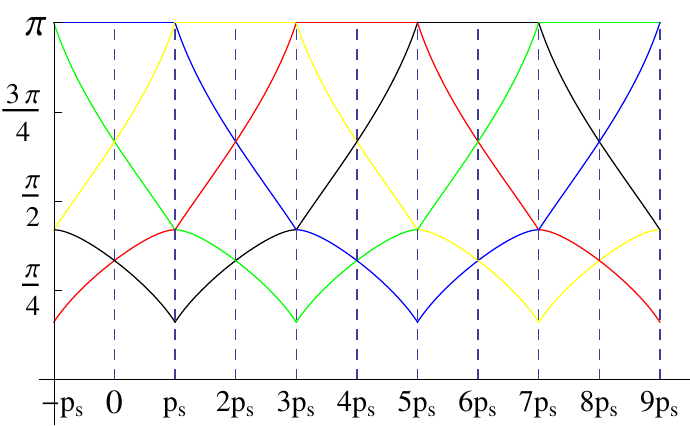}
\label{subfig:psiv}
}
\subfigure[]
{\includegraphics[width=7cm]{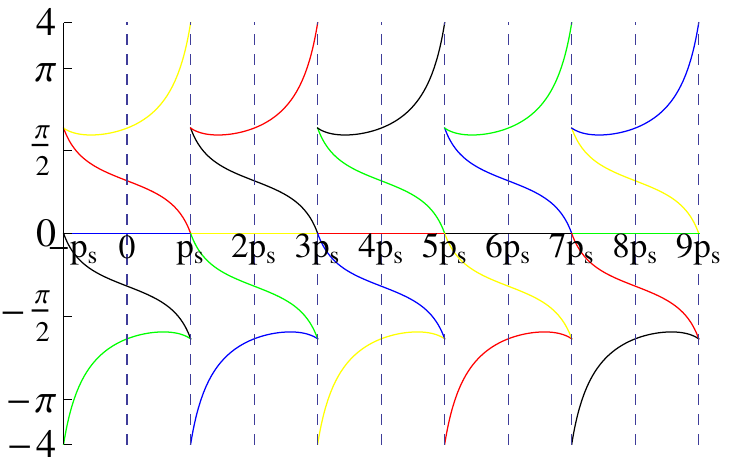}
\label{subfig:dpsiv}
}
\caption{The relative angles $\psi_i$ (left) and $\dot{\psi_i}$ (right) along the loop of the boundary of the convex pentagons in the clockwise direction starting at the red point located at $ p(-p_s)$ in Figure~\ref{subfig:pentcurve} where $p$ is from \eqref{eq:pcur}.
%
%
The colouring of the curves $i=1,\dots, 5$ is green, red, blue, black, yellow (standard colour code). 
The dashed vertical lines denotes the boundary of each fundamental piece, 
and the finite jumps for $\dot \psi_i$ that occur at every second boundary are the result of the corner points from Figure~\ref{subfig:pentcurve}.
}\label{fig:psivex}
\end{figure}

The derivatives $\dot{\psi_i}$ along the boundary of the convex region have finite jumps at $kp_s$ when $k$ is odd, these jumps are shown in 
Figure~\ref{subfig:dpsiv} and corresponds to the corners of the white loop from Figure \ref{subfig:pentcurve}. The corresponding shapes at these corners are the isosceles triangles with side lengths $\{1,2,2\}$. For even $k$ in $kp_s$, the shapes are the trapeziums with side lengths $\{1,1,1,2\}$ and the transition from one fundamental region to another at this point is smooth, as seen in Figure~\ref{subfig:dpsiv}.
By integrating $\dot{\theta}$ in \eqref{eq:thetadot} we obtain a continuous function $\theta(t)$ (though not smooth at $kp_s$ for odd $k$), 
and the overall rotation generated by this loop is found to be 
$\Delta\theta\approx 0.33117\approx 19.01^{\circ}$. 
This is the maximal overall rotation achievable after one period if we restrict the allowed shapes to the space of all equilateral convex pentagons.

\section{Conclusion}
The optimal way for the equilateral pentagon at zero total angular momentum to achieve maximal overall rotation after one period of 
a periodic contractible loop is to follow the $B=0$ contour around the regular pentagram. This results in an overall rotation of $\Delta\theta\approx 0.78837$ radians or $\approx 45.17^\circ$.

If we restrict to the subset of simple equilateral pentagons, then the $B=0$ contour around the regular convex pentagon provides the maximal overall rotation with $\Delta\theta\approx 0.49147$ radians or $\approx 28.16^\circ$, this is the black loop found in Figure~\ref{subfig:pentcurve}. Finally, if we further restrict to the subset of equilateral convex pentagons, then the maximal overall rotation reduces to $\Delta\theta\approx 0.33117$ radians or $\approx 19.01^\circ$. 
The fact that this loop is smaller and entirely contained within the larger $B=0$ loop around the regular convex pentagon (see the white loop in Figure~\ref{subfig:pentcurve}) explains why the overall rotation is smaller for this loop.

An intuitive explanation of why the loop around the pentagram gives a largest value is that 
the moment of inertia for the pentagram is the global minimum. Hence the magnetic 
field $B$ tends to be bigger there, as compared to the magnetic field near the regular convex pentagon,
which has maximal moment of inertia. 

\newpage
\section*{Appendix}
Movies of .mp4 format for optimal shape changes are available at\\\\
\url{http://www.maths.usyd.edu.au/u/williamt/pentagon/Beq0_pentagram.mp4}\\\\
for the optimal loop around the pentagram,\\\\
\url{http://www.maths.usyd.edu.au/u/williamt/pentagon/Beq0_convex_pentagon.mp4}\\\\
for the optimal loop around the pentagon, and\\\\
\url{http://www.maths.usyd.edu.au/u/williamt/pentagon/subset_convex_pentagon.mp4}\\\\
for the optimal loop for convex pentagons.\\\\
\newline
Interactive Mathematica CDF files are available at\\\\
\url{http://www.maths.usyd.edu.au/u/williamt/pentagon/manipulate.cdf}\\\\
for finding the shape corresponding to a point in the fundamental region, and\\\\
\url{http://www.maths.usyd.edu.au/u/williamt/pentagon/animations.cdf}\\\\
for animations of the optimal shape changes.\\\\
\newline
Download the Wolfram CDF player at\\\\
\url{http://www.wolfram.com/cdf-player/}

\section*{Acknowledgement}
This research was supported by the ARC grant LP100200245.
We would like to thank Leon Poladian for useful comments.
HRD would like to thank the Department of Applied Mathematics
of the University of Colorado at Boulder for their hospitality.
\newpage
\bibliographystyle{siam}
\bibliography{pentagon}

\begin{thebibliography}{10}

\bibitem{Bloch03}
{\sc A.~M. Bloch}, {\em Nonholonomic mechanics and control}, vol.~24 of
  Interdisciplinary Applied Mathematics, Springer-Verlag, New York, 2003.

\bibitem{genus4a}
{\sc R.~Curtis and M.~Steiner}, {\em Configuration spaces of planar pentagons},
  American Mathematical Monthly, 114 (2007), pp.~183--201(19).

\bibitem{genus4b}
{\sc T.~F. Havel}, {\em Some examples of the use of distances as coordinates
  for euclidean geometry}, Journal of Symbolic Computation, 11 (1991),
  pp.~579--593.

\bibitem{MacKay03}
{\sc T.~J. Hunt and R.~S. MacKay}, {\em Anosov parameter values for the triple
  linkage and a physical system with a uniformly chaotic attractor},
  Nonlinearity, 16 (2003).

\bibitem{Hyde97}
{\sc S.~Hyde}, {\em The Language of Shape}, Elsevier, 1997.

\bibitem{KaneScher69}
{\sc T.~Kane and M.~Scher}, {\em A dynamical explanation of the falling cat
  phenomenon}, International Journal of Solids {\&} Structures, 5 (1969),
  pp.~663--670.

\bibitem{LittlejohnReinsch97}
{\sc R.~G. Littlejohn and M.~Reinsch}, {\em Gauge fields in the separation of
  rotations and internal motions in the n-body problem}, Rev. Mod. Phys., 69
  (1997), pp.~213--275.

\bibitem{Marsden98symmetriesin}
{\sc J.~Marsden}, in Motion, Control, and Geometry: Proceedings of a Symposium,
  {Board on Mathematical Sciences} and {National Research Council}, eds.

\bibitem{MMR90}
{\sc J.~Marsden, R.~Montgomery, and T.~Ratiu}, {\em Reduction, symmetry, and
  phases in mechanics}, Mem. Amer. Math. Soc., 88 (1990), pp.~iv+110.

\bibitem{Montgomery93cat}
{\sc R.~Montgomery}, {\em Gauge theory of the falling cat}, in Dynamics and
  control of mechanical systems ({W}aterloo, {ON}, 1992), vol.~1 of Fields
  Inst. Commun., Amer. Math. Soc., Providence, RI, 1993, pp.~193--218.

\bibitem{PuttermanRaz08}
{\sc E.~Putterman and O.~Raz}, {\em The square cat}, American Journal of
  Physics, 76 (2008), pp.~1040--1044.

\bibitem{YangKrish}
{\sc R.~Yang and P.~S. Krishnaprasad}, {\em On the geometry and dynamics of
  floating four-bar linkages}, Dynamics and Stability of Systems, 9 (1994),
  pp.~19--45.

\end{thebibliography}
\end{document}